\documentclass[12pt]{amsart}

\author{Paul \textsc{Poncet}}
\address{CMAP, \'{E}cole Polytechnique, Route de Saclay, 91128 Palaiseau Cedex, France}

\email{poncet@cmap.polytechnique.fr}

\usepackage[colorinlistoftodos,bordercolor=orange,backgroundcolor=orange!20,linecolor=orange,textsize=scriptsize]{todonotes}
\usepackage{stmaryrd}
\usepackage{amssymb}
\usepackage{amsmath}
\usepackage{yhmath}
\usepackage{calrsfs} % nouvelles lettres calligraphiques (commande \mathrsfs)
\usepackage{times} % sans que je comprenne pourquoi, les recherches textuelles dans les documents pdf fonctionnent avec cette fonte
\usepackage{tikz}

\def\twoheaddownarrow{\rlap{$\downarrow$}\raise-.5ex\hbox{$\downarrow$}}%%
\def\twoheaduparrow{\rlap{$\uparrow$}\raise.5ex\hbox{$\uparrow$}}%%
\newcommand{\A}{\mathsf{A}}
\newcommand{\B}{\mathsf{B}}
\newcommand{\C}{\mathsf{C}}
\newcommand{\D}{\mathsf{D}}

\newcommand{\F}{\mathsf{F}}
\newcommand{\I}{\mathsf{I}_{\mathsf{Z}}} %{\downarrow\!\!\mathsf{z}}
\newcommand{\Z}{\mathsf{Z}}

\newtheorem*{theorem*}{Theorem}

\newtheorem{theorem}{Theorem}[section]
\newtheorem{corollary}[theorem]{Corollary}
\newtheorem{proposition}[theorem]{Proposition}
\newtheorem{lemma}[theorem]{Lemma}

\theoremstyle{definition}
\newtheorem{definition}[theorem]{Definition}
\newtheorem*{definition*}{Definition}

\newtheorem{remark}[theorem]{Remark}
\newtheorem{example}[theorem]{Example}
\newtheorem{problem}[theorem]{Problem}
\newenvironment{acknowledgements}[1][]{\par\noindent\textbf{Acknowledgements#1.} }{\par}

\begin{document}

\title{Transporting continuity properties \\ from a poset to its subposets}
\date{\today}

 %Continuous lattices and posets, applications 


\begin{abstract}
We identify two key conditions that a subset $A$ of a poset $P$ may satisfy to guarantee the transfer of continuity properties from $P$ to $A$. 
We then highlight practical cases where these key conditions are fulfilled. 
Along the way we are led to consider subsets of a given poset $P$ whose way-below relation is the restriction of the way-below relation of $P$, which we call way-below preserving subposets. 
As an application, we show that every conditionally complete poset with the interpolation property contains a largest continuous way-below preserving subposet. 
Most of our results are expressed in the general setting of $\Z$ theory, where $\Z$ is a subset system. 
\end{abstract}

\maketitle

%%%%%%%%%%%%%%%%%%%%%%
%%%%%%%%%%%%%%%%%%%%%%
%%%%%%%%%%%%%%%%%%%%%%
%%%%%%%%%%%%%%%%%%%%%%
\section{Introduction}

In the article \cite{Venugopalan86} Venugopalan showed that, if $A$ is a subset of a conditionally directed-complete poset $P$ such that the inclusion map $A \to P$ is Scott-continuous and the right part of a Galois connection, then $A$ is a continuous poset if $P$ is itself continuous. 
%This is a powerful way of transferring continuity properties from a poset to one of its subsets. 
In this paper, our goal is to provide a series of results of the same flavor. 
For this purpose, we identify two key properties: preservation of directed sups, and \textit{refinement}. 
Preserving directed sups basically means for a map to be Scott-continuous; for a subset that the inclusion map is Scott-continuous. 

The refinement property has a slightly more technical definition. 
We say that an order-preserving map $f : P \to P'$ between posets has the \textit{refinement property} if, whenever $x \in P$ and $f(x) \leqslant \bigvee_{P'} D'$ for some directed subset $D'$ of $P'$ with sup, there is a directed subset $D$ of $P$ with sup such that ${ f(D) } \subseteq { \downarrow_{P'}\!\! D' }$ and $x \leqslant \bigvee_P D$. 
A subset $A$ of a poset $P$ has the \textit{refinement property} if the inclusion map $A \to P$ has the refinement property. 

Combining the preservation of directed sups and the refinement property 
guarantees the transfer of continuity properties from $P$ to $A$, as stated by our main theorem.

\begin{theorem*}[Main Theorem, shortened version]\label{thm:mainshort}
Let $A$ be a subset of a poset $P$ that preserves directed sups and has the refinement property. 
Then $A$ is a way-below preserving subposet of $P$ and the Scott-open subsets of $A$ coincide with the subsets induced by the Scott-open subsets of $P$. 
Moreover, if $P$ is continuous, then $A$ is continuous.
\end{theorem*}

Here by a \textit{way-below preserving subposet} of a poset $P$ we mean a subset $A$ preserving directed sups 
and such that the way-below relation $\ll_A$ on $A$ coincides with the restriction to $A \times A$ of the way-below relation $\ll_P$ on $P$. 

%For this purpose, we shall identify two key properties that a subset $A$ of a poset $P$ may satisfy to guarantee the transfer of continuity properties from $P$ to $A$, in the spirit of Venugopalan's result. We shall then investigate practical cases where these key properties are fulfilled. 

%Recognizing a poset as being continuous is a valuable information, as the area of domain theory demonstrates. 
Our motivation comes from the observation that ad hoc technical proofs are often required that a given subset $A$ at stake be a continuous poset, using various hypotheses on $A$ in addition to continuity properties of the ambient poset $P$. 
Such situations can be seen in the domain theoretic literature; for instance Lawson and Xu \cite{Lawson04} investigated continuity properties of principal filters, principals ideals, and order intervals of a poset. 
This can also be the case in mathematical areas where domain theory has applications. 
For instance Keimel \cite{Keimel09} considered continuity properties of cones of a finite-dimensional vector space $\mathbb{R}^n$. 
Recall that a \textit{cone} $C$ of $\mathbb{R}^n$ is a subset such that $C + C \subseteq C$, 
$\lambda C \subseteq C$ for all $\lambda \geqslant 0$, and $C \cap (-C) = \{ 0 \}$. 
It induces a partial order $\leqslant_C$ on $\mathbb{R}^n$ defined by $x \leqslant_C y$ if $y - x \in C$. 
Keimel asserted without proof that $\mathbb{R}^n$ equipped with this partial order becomes a continuous and dually continuous poset if $C$ has nonempty interior. 
In this case, he derived that other subsets of $\mathbb{R}^n$ are continuous posets too, including $C$ and its topological interior. 

Another situation is that of idempotent analysis, where methods from domain theory have proved their usefulness at various occasions, see e.g.\ Akian \cite{Akian99}, Akian and Singer \cite{Akian03}, Lawson \cite{Lawson04b}, Poncet \cite{Poncet11, Poncet12b, Poncet12c}. 
An algebraic structure at stake in this area is that of \textit{idempotent semifield}, that is a semiring $(K, +, 0, \times, 1)$ in which every nonzero element has a multiplicative inverse and such that $x + x = x$ for all $x$; it gives rise to a partial order $\leqslant_K$ defined by $x \leqslant_K y$ if $x + y = y$. 
In \cite[Chapter~3]{Poncet11}, we showed how continuity properties of (subsets of) $K$ come into play when considering Riesz like idempotent representation theorems of ``linear" forms defined on a $K$-module. 

A class of algebraic structures that offers a potential framework for unifying classical and idempotent analysis is that of \textit{inverse semigroups}, which form a well studied class of semigroups endowed with a natural partial order. 
%Since continuity properties have proved to play an important role in idempotent analysis, 
It was a source of motivation of our work \cite{Poncet12a} to study continuity properties of inverse semigroups; we proved that, under mild conditions on an inverse semigroup $S$, the latter is continuous as a poset if and only if its subsemigroup made of idempotent elements (which is a semilattice) is itself continuous. 

%Dire aussi peut-etre que de nouvelles applications des continuous posets apparaissent : Furber (C* algebras), Keimel, Martin, Poncet... ou de tels lemmes techniques sont justement indispensables.

%In all the previous situations, technical lemmata are required to transfer continuity properties from one poset to another. 
%It happens that all these different situations can be managed in a similar fashion; this is the main goal of this paper to highlight this and enable users of domain theory to deduce continuity properties of posets and their subsets more easily. 
%We show that our main theorem is applicable in a wide variety of situations. 
We show that our main theorem is widely applicable; especially, the hypotheses of the theorem are met in the following cases: 
\begin{itemize}
  \item if $A$ is a lower set in a continuous semilattice;
  \item if $A$ is a lower set in a continuous conditionally directed-complete poset;
  \item if $A$ is a Scott-convex subset (defined as the intersection of a lower set with a Scott-open set) in a meet-continuous semilattice;
  \item if $A$ is a Scott-open subset of a poset;
  \item if $A$ is the image of a conditionally directed-complete poset by a Scott-continuous projection with the refinement property;
  \item if $A$ is a kernel retract of a poset induced by a Scott-continuous kernel retraction;
  \item if $A$ is already a continuous way-below preserving subposet of a poset.
\end{itemize}

We apply our main theorem and especially the case where $A$ is a kernel retract to a question that was another source of motivation at the start of this work and relates to the interpolation property of posets. % question of the existence of continuous subposets within a poset with the interpolation property. 
%Indeed, another source of motivation at the start of this work has been to study the interpolation property of posets. 
This property is a necessary condition for a poset to be continuous, 
and is also crucial for deriving many important results of the theory of continuous posets; 
however, little has been done on (not necessarily continuous) posets with the interpolation property. 
A notable exception is the recent work of Mao and Xu \cite{Mao17}: these authors showed that 
a conditionally complete poset with the interpolation property is continuous 
if and only if $\twoheaddownarrow x = \twoheaddownarrow y$ implies $x = y$ for all $x, y$. 
With this work in mind, we give examples of posets with the interpolation property that are not continuous, 
and we examine whether a poset with the interpolation property 
can contain subposets that are continuous posets with respect to the induced partial order. 
This leads us to the following result: 

\begin{theorem*}
Every conditionally complete poset with the interpolation property contains a largest continuous way-below preserving subposet, 
which is itself conditionally complete. 
\end{theorem*}

Given such a conditionally complete poset $P$ with the interpolation property, 
the strategy of proof is to apply our main theorem to the subset $A$ defined as the kernel retract of $P$ 
induced by the Scott-continuous kernel operator $k : P \to P, x \mapsto \bigvee \twoheaddownarrow x$. 

While this introductory section has been expressed in the language of classical domain theory for the sake of simplicity, 
we actually adopt in this paper the general $\Z$ framework of Wright, Wagner, and Thatcher \cite{Wright78}. 
Recall that it mainly consists in replacing directed subsets used in the definitions by some general family of subsets $\Z(P)$, for each poset $P$. 
This choice is not only guided by a sake of generality: we need it for some applications cited above, notably for Riesz like representation theorems in idempotent analysis. 
Another reason is that, in our results, several families of subsets are sometimes intertwined, hence the $\Z$ framework provides ease of formulation and clarity. 

The paper is organized as follows. 
In Section~\ref{sec:interp} we introduce partially ordered sets (posets) and their continuity, when accompanied with a subset system $\Z$; 
we recall concepts such as $\Z$-subsets, $\Z$-interpolation property, $\Z$-continuity, $\Z$-basis, $\Z$-sup-preserving maps, and Galois connections. 
We also give examples of posets with the interpolation property that are not continuous. 
In Section~\ref{sec:maps} we focus on order-preserving maps between posets; we introduce the property of $\Z$-refinement and make the link with the more usual property of $\Z$-below preservation. 
In Section~\ref{sec:main} we specialize the definitions of the previous section to subsets of posets seen as inclusion maps, 
and we notably formulate and prove our Main Theorem in its full generality. 
In passing, we examine the question of transferring the continuity properties of a covering family of subposets of a poset to the poset itself. 
In the next five sections, we apply our Main Theorem to the various situations evoked above. 
In particular, we consider the case of lower subsets in Section~\ref{sec:lower}, 
order-convex subsets and $\Z$-Scott-convex subsets in Section~\ref{sec:convex}, 
upper subsets and Scott open subsets in Section~\ref{sec:upper}, 
direct images of order-preserving maps in Section~\ref{sec:left}, 
kernel retracts in Section~\ref{sec:kernels}. 
We conclude this article with Section~\ref{sec:largest}, where we generalize a theorem by Mao and Xu \cite{Mao17}
and show in particular that every conditionally complete poset with the interpolation property contains a largest (conditionally complete) continuous way-below preserving subposet.

%%%%%%%%%%%%%%%%%%%%%%
%%%%%%%%%%%%%%%%%%%%%%
%%%%%%%%%%%%%%%%%%%%%%
%%%%%%%%%%%%%%%%%%%%%%
\section{Continuous posets and the interpolation property}\label{sec:interp}

\subsection{Posets}

A \textit{partially ordered set} (\textit{poset}, for short) is a set $P$ equipped with a partial order $\leqslant_P$, i.e.\ an antisymmetric, reflexive, transitive binary relation. 
Let $A \subseteq P$ and $x \in P$. 
We write $\downarrow_{P}\!\! A$ for the subset $\{ y \in P : y \leqslant_P a \mbox{ for some $a \in A$ } \}$, and $\downarrow_{P}\!\! x$ as a shorthand for the \textit{principal ideal} $\downarrow_{P}\!\! \{ x \}$. 
We write $\uparrow_{P} \!\! A$ for the subset $\{ y \in P : a \leqslant_P y \mbox{ for some $a \in A$ } \}$, and $\uparrow_{P}\!\! x$ for the \textit{principal filter} $\uparrow_{P}\!\!\{ x \}$. 
The subset $A$ is a \textit{lower set} if $A = { \downarrow_{P}\!\! A }$, an \textit{upper set} if $A = { \uparrow_{P}\!\! A }$, and is \textit{order-convex} if $[x, z]_P \subseteq { A }$ for all $x, z \in A$, where $[x, z]_P$ denotes $\{ y \in P : x \leqslant_P y \leqslant_P z \}$. 
Equivalently, $A$ is order-convex if and only if it is the intersection of a lower subset with an upper subset; in particular, lower subsets and upper subsets are order-convex. 
An \textit{upper bound} of $A$ is an element $u \in P$ such that $a \leqslant_P u$ for all $a \in A$. 
We write $A^{\uparrow}$ for the subset of upper bounds of $A$. 
The \textit{supremum} (or \textit{sup}) of $A$, if it exists, is the least upper bound of $A$, i.e.\ it is an upper bound $u_0$ of $A$ such that $u_0 \leqslant_P u$ for every upper bound $u$ of $A$; we denote it by $\bigvee_P A$. 
The notions of \textit{lower bound} and \textit{infimum} (or \textit{inf}) are defined dually; we denote by $A^{\downarrow}$ the subset of lower bounds of $A$, and by $\bigwedge_P A$ the inf of $A$ in $P$, if it exists. 
When the context is clear, we remove the subscript $P$ and write $\leqslant$, $\downarrow\!\! A$, $\downarrow\!\! x$, $\bigvee A$, etc. instead of $\leqslant_P$, $\downarrow_{P}\!\! A$, $\downarrow_{P}\!\! x$, $\bigvee_P A$, etc. 

Given two posets $P$ and $Q$, a map $f : P \to Q$ is \textit{order-preserving} if $f(x) \leqslant_Q f(y)$ whenever $x \leqslant_P y$, or equivalently if ${ f(\downarrow_{P}\!\! A) } \subseteq { \downarrow_{Q}\!\! f(A) }$ for all $A \subseteq P$.
We denote by $f^{\circ}$ the order-preserving map $P \to f(P)$, $x \mapsto f(x)$, called the \textit{corestriction of $f$}. 

An order-preserving map $f : P \to Q$ is \textit{quasi-invertible} if $f \circ g \circ f^{\circ} = f$, for some order-preserving map $g : f(P) \to P$; the map $g$ is called a \textit{quasi-inverse} of $f$ and is necessarily injective (it is even an \textit{order-embedding}, see the definition at the beginning of Section~\ref{sec:maps}). 
A \textit{projection} on $P$ is an order-preserving map $p : P \to P$ such that $p \circ p = p$. 
For instance, if $f : P \to Q$ is quasi-invertible with quasi-inverse $g$, then $g \circ f^{\circ}$ is a projection on $P$ and $f^{\circ} \circ g$ is a projection on $f(P)$.  

\subsection{Subset systems}

A \textit{subset system} is a function $\Z$ that assigns to each poset $P$ a collection $\Z(P)$ of subsets of $P$ such that 
\begin{enumerate}
	\item[$i)$] there exists a poset $P$ such that $\Z(P)$ has a nonempty element; 
	\item[$ii)$] $f(Z) \in \Z(Q)$, for every order-preserving map $f : P \to Q$ and every $Z \in \Z(P)$. 
%	\item[$iii)$] there exists a poset $P$ such that $\Z(P)$ contains a set with cardinality greater than one. 
\end{enumerate}
The elements of $\Z(P)$ are called the \textit{$\Z$-subsets} of $P$. 
We write $\Z^{*}$ for the subset system defined by $\Z^{*}(P) = \Z(P) \setminus \{ \emptyset \}$ for every poset $P$. 

The two conditions $i)$ and $ii)$ correspond to the original definition of a subset system given by Wright et al.\ \cite{Wright78}. Taken together, they ensure that the singletons of $P$ are $\Z$-subsets of $P$, for every poset $P$. 
Condition $ii)$ means that $\Z$ is a covariant functor from the category of posets to the category of sets with $\Z(f) : \Z(P) \to \Z(Q)$ defined by $\Z(f)(Z) = f(Z)$ where $Z \in \Z(P)$, for every order-preserving map $f : P \to Q$. 

We say that a subset system $\Z$ is \textit{proper} if it satisfies the following additional condition:
\begin{enumerate}
	\item[$iii)$] there exists a poset $P$ such that $\Z(P)$ contains a set with at least two elements. 
\end{enumerate}
Condition $iii)$, which obviously implies Condition $i)$, was added later by Banaschewski and Nelson \cite{Banaschewski82} and Baranga \cite{Baranga96}; it is quite desirable, as $ii)$ and $iii)$ taken together imply that every pair $\{ x, y \}$ such that $x \leqslant y$ is a $\Z$-subset \cite[Theorem~1.6]{Baranga96}, and consequently that every $\Z$-sup-preserving map is order-preserving \cite[Corollary~1.7]{Baranga96} (see Paragraph~\ref{subsec:sl} below for the definition of a $\Z$-sup-preserving map). 

%To this definition, first given by Wright et al.\ \cite{Wright78}, we add a third (unusual but useful in the framework of this paper) condition: 
%\begin{enumerate}
%	\item[$iii)$] the empty set is not in $\Z(P)$, for all posets $P$. 
%\end{enumerate}

The suggestion of \cite{Wright78} to apply subset systems to the theory of continuous posets was followed by Nelson \cite{Nelson81}, Banaschewski and Nelson \cite{Banaschewski82}, Novak \cite{Novak82b}, Bandelt \cite{Bandelt82}, Bandelt and Ern\'e \cite{Bandelt83}, \cite{Bandelt84}, and this research was carried on by Venugopalan \cite{Venugopalan86}, \cite{Venugopalan88}, Xu \cite{Xu95}, Baranga \cite{Baranga96}, Menon \cite{Menon96}, %Shi and Wang \cite{Shi96}, 
Ern\'e \cite{Erne99}, \cite{Erne01}, Zhao \cite{Zhao15}, Yuan and Li \cite{Yuan19}, among others. 

A nonempty subset $D$ of a poset is \textit{directed} if, for all $x, y \in D$, one can find $d \in D$ such that $x \leqslant d$ and $y \leqslant d$. The proper subset system $\D$ that selects the directed subsets of each poset is behind the classical theory of continuous posets and domains, see the monographs by Gierz et al.\ \cite{Gierz03} and Goubault-Larrecq \cite{Goubault13}. 
Here are some further examples of proper subset systems: 
\begin{itemize}
%	\item For $\Z(P)$ the set of Frink ideals (resp. Cauchy ideals) of $P$, see
%	\item If $\Z(P)$ is the set of countably generated ideals, one gets the object studied in \cite{Banaschewski81}.  (In this paper $P$ is lattice with top and bottom such that every countable set has a join.)
  \item $\A$ (resp.\ $\A^{*}$) selects all (resp.\ nonempty) subsets. It works well for investigating completely distributive lattices, see Raney \cite{Raney52}, \cite{Raney53}, Ern\'e et al.\ \cite{Erne06}. 
%	\item The case where $\Z(P)$ is the set of filtered subsets of $P$ was used for instance by G.\ Gerritse \cite{Gerritse97}, Jonasson \cite{Jonasson98}, Akian and Singer \cite{Akian03}. See also \cite{Poncet12b}. 
  \item $\B$ selects upper-bounded subsets. 
  \item $\C$ selects \textit{chains}, i.e.\ subsets $C$ such that $x \leqslant y$ or $y \leqslant x$ whenever $x, y \in C$. See Markowsky and Rosen \cite{Markowsky76b}, and Markowsky \cite{Markowsky76}, \cite{Markowsky77}, \cite{Markowsky81a}, \cite{Markowsky81b}. %Using the Hausdorff maximality theorem, relations between directed subsets and chains were explored by Iwamura \cite{Iwamura44}, Bruns \cite{Bruns67}, and Markowsky \cite{Markowsky76}. 
  See also Ern\'e \cite[p.\ 54]{Erne99}. 
%    \item $\E^{*}$ selects singletons. % the way-above relation $y \gg x$ (defined below) reduces to the partial order $y \geqslant x$.} %And yet it will play a  crucial in Chapter~\ref{ch:maxmes}. 
  \item $\F$ selects finite subsets, see Martinez \cite{Martinez72}, %, and is linked with the abstract convexity theory developed in van de Vel's monograph \cite{vanDeVel93}. 
	Frink \cite{Frink54}, Ern\'e \cite{Erne81}. 
\end{itemize}

Given a subset system $\Z$, we have the inclusion 
\[
\Z(A) \subseteq { \{ Z \in \Z(P) : Z \subseteq A \} }
\]
for every subset $A$ of a poset $P$, and equality holds for $\Z \in \{ \A, \A^*, \C, \D, \F \}$. 
For arbitrary $\Z$, equality can also be reached in special cases, as the following result testifies. 

\begin{proposition}\label{prop:stable}
Let $\Z$ be a subset system, $P$ be a poset, $x \in P$, $Z$ be a $\Z$-subset of $P$, and $A \subseteq P$. 
The following assertions hold: 
\begin{enumerate}
  \item\label{prop:stable1} If $x \in Z$, then $Z \cap \uparrow\!\! x$ is a $\Z$-subset of $\uparrow\!\! x$ (hence a $\Z$-subset of $P$). 
  \item\label{prop:stable2} If $x \in Z$, then $Z \cap \downarrow\!\! x$ is a $\Z$-subset of $\downarrow\!\! x$ (hence a $\Z$-subset of $P$). 
  \item\label{prop:stable3} If $x \in Z^{\downarrow}$, then $Z$ is a $\Z$-subset of $\uparrow\!\! x$. 
  \item\label{prop:stable4} If $x \in Z^{\uparrow}$, then $Z$ is a $\Z$-subset of $\downarrow\!\! x$. 
  \item\label{prop:stable5} If $Z^{\downarrow} \cap A \neq \emptyset$ and $A$ is an upper set, then $Z$ is a $\Z$-subset of $A$. 
  \item\label{prop:stable6} If $Z^{\uparrow} \cap A \neq \emptyset$ and $A$ is a lower set, then $Z$ is a $\Z$-subset of $A$. 
  \item\label{prop:stable7} If $Z \subseteq { A }$ and $A = { f(P') }$ for some quasi-invertible map $f : P' \to P$, then $Z$ is a $\Z$-subset of $A$. 
\end{enumerate}
\end{proposition}

\begin{proof}
Cases \eqref{prop:stable2}, \eqref{prop:stable4}, and \eqref{prop:stable6} are dual versions of, respectively, \eqref{prop:stable1}, \eqref{prop:stable3}, and \eqref{prop:stable5}, so we only prove \eqref{prop:stable1}, \eqref{prop:stable3}, \eqref{prop:stable5}, and \eqref{prop:stable7}. 
Let $f_x : P \to { \uparrow\!\! x }$ be the order-preserving map defined by $f_x(y) = y$ if $x \leqslant y$, and $f_x(y) = x$ otherwise. 

Case~\eqref{prop:stable1}: 
If $x \in Z$, then $f_x(Z) = Z \cap \uparrow\!\! x$, so $Z \cap \uparrow\!\! x$ is a $\Z$-subset of $\uparrow\!\! x$. 

Case~\eqref{prop:stable3}: 
If $x \in Z^{\downarrow}$, then $f_x(Z) = Z$, so $Z$ is a $\Z$-subset of $\uparrow\!\! x$. 

Case~\eqref{prop:stable5}: 
Let $A$ be an upper subset of $P$ such that $Z^{\downarrow} \cap A$ is nonempty. 
Let $y \in { Z^{\downarrow} \cap A }$. 
Then $Z$ is a $\Z$-subset of $\uparrow\!\! y$ by \eqref{prop:stable3}. 
Moreover, $\uparrow\!\! y \subseteq { A }$ since $A$ is an upper set, so $Z$ is a $\Z$-subset of $A$. 

Case~\eqref{prop:stable7}: 
Let $g : A \to P'$ be a quasi-inverse of $f$, and let $p : A \to A$, $x \mapsto f(g(x))$. 
If $z \in Z  \subseteq { A } = { f(P') }$, then $z = f(x')$ for some $x' \in P'$, so $p(z) = f(g(f(x'))) = f(x') = z$. 
This shows that $p(Z) = Z$. 
Since $p(Z)$ is a $\Z$-subset of $A$, $Z$ is a $\Z$-subset of $A$. 
\end{proof}

%\begin{definition}
%A subset system $\Z$ is \textit{intrinsically defined} if
%\[
%\Z(A) = \{ Z \in \Z(P) : Z \subseteq A \}, 
%\]
%for every poset $P$ and every subset $A$ of $P$ equipped with the induced partial order. 
%\end{definition}

%All subset systems listed above but $\B$ are intrinsically defined. 

In the remaining part of this paper, $\Z$ will always denote a proper subset system.

\subsection{Continuity and the interpolation property}

Let $P$ be a poset. 
Then $P$ is said to be \textit{$\Z$-complete} if every $\Z$-subset of $P$ has a sup. 
An element $x \in P$ is \textit{$\Z$-below} $y \in P$, written $x \ll^{\Z}_P y$, if, for every $\Z$-subset $Z$ with sup, $y \leqslant \bigvee Z$ implies $x \in { \downarrow\!\! Z }$. 
We write $\twoheaddownarrow^{\Z}_P x$ for the subset $\{ y \in P : y \ll^{\Z}_P x \}$. %, and $\twoheaduparrow{}^{\Z}_P x$ for the subset $\{ y \in P : x \ll^{\Z}_P y \}$. 
When the context is clear, we remove the subscript $P$ and write $\ll^{\Z}$, $\twoheaddownarrow^{\Z} x$ instead of $\ll^{\Z}_P$, $\twoheaddownarrow^{\Z}_P x$. 

An element $x \in P$ is \textit{$\Z$-compact} if $x \ll^{\Z} x$. 
The poset $P$ is \textit{$\Z$-continuous} if, for each $x \in P$, the subset $\twoheaddownarrow^{\Z} x$ contains a $\Z$-subset of $P$ whose sup is $x$. 
We emphasize that $\Z$-completeness is \textit{not} part of our definition of $\Z$-continuity, in agreement with Baranga \cite{Baranga96}, but unlike Bandelt and Ern\'e \cite{Bandelt83} or Venugopalan \cite{Venugopalan86}. 
The notion of $\Z$-continuity in the sense of this paper coincides with \textit{$\Z_{\vee}$-precontinuity} in Ern\'e \cite{Erne99}. 

In the literature (see e.g.\ \cite{Gierz03}), the $\D$-below relation is known as the \textit{way-below} relation, and $\D$-continuous posets are known as \textit{continuous posets}; yet, we shall stick to the $\D$ designation in the frame of this paper. 

The following example is taken from \cite[Example~2.7]{Venugopalan86}. 

\begin{example}[Venugopalan]
Let $\mathbb{N}$ be the set of natural numbers ordered by divisibility. 
Then $\mathbb{N}$ is an $\F$-continuous poset, and $m \ll^{\F} n$ if and only if $m$ is a prime power divisor of $n$. 
\end{example}

The following well-known result will be used repeatedly and often implicitly along the paper. 

\begin{lemma}
Let $P$ be a poset, and let $w, x, y, z \in P$. 
\begin{enumerate}
  \item if $x \ll^{\Z} y$, then $x \leqslant y$;
  \item if $w \leqslant x \ll^{\Z} y \leqslant z$, then $w \ll^{\Z} z$. 
\end{enumerate}
\end{lemma}

\begin{proof}
See e.g.\ \cite[Proposition~2.2]{Baranga96}. 
\end{proof}

A poset $P$ has the \textit{$\Z$-interpolation property}, or is \textit{$\Z$-interpolating} if, for all $x, z \in P$ with $x \ll^{\Z} z$, there exists some $y \in P$ such that $x \ll^{\Z} y \ll^{\Z} z$. 
A $\Z$-continuous poset with the $\Z$-interpolation property is called \textit{strongly $\Z$-continuous} (\cite{Novak82b}, \cite{Bandelt83}). 

In a poset, a \textit{$\Z$-ideal} is a subset $I$ such that $I = { \downarrow\!\! Z }$ for some $\Z$-subset $Z$. 
Note that a poset $P$ is $\Z$-continuous if and only if, for each $x \in P$, the subset $\twoheaddownarrow^{\Z} x$ is a $\Z$-ideal and has a sup equal to $x$ (see e.g.\ \cite[Theorem~2.6]{Baranga96}). 
We denote by $\I(P)$ the poset of $\Z$-ideals of $P$ ordered by inclusion. 
%The subset system made of the $\Z$-ideals is denoted by $\I$. 
The proper subset system $\Z$ is \textit{union-complete} if, for every $\mathrsfs{V} \in \Z(\I(P))$, $\bigcup \mathrsfs{V} \in \I(P)$. 
All subset systems mentioned above but $\C$ (see \cite{Banaschewski82}) are union-complete. 
% this encodes the fact that finite unions of finite sets are finite, directed unions of directed subsets are directed, etc. 
%While it remains an open problem to exhibit a $\Z$-continuous poset with respect to some subset system $\Z$ that does not satisfy the $\Z$-interpolation property, the following result gives a sufficient condition on $\Z$

\begin{remark}\label{rk:unionc}
A common way of taking benefit from union-completeness is the following, when one disposes of a family $(I_x)_{x \in P}$ of $\Z$-ideals of a poset $Q$ indexed by a poset $P$, and such that ${ x \leqslant y } \Rightarrow { I_x \subseteq I_y }$, for all $x, y \in P$. 
The map $i : P \to \I(Q)$, $x \mapsto I_x$ is order-preserving, hence the direct image $i(Z)$ of a $\Z$-subset $Z$ of $P$ by $i$ is a $\Z$-subset of $\I(Q)$, i.e.\ $i(Z) \in \Z(\I(Q))$. 
Since ${ \bigcup i(Z) } = { \bigcup_{x \in Z} I_x }$, union-completeness of $\Z$ gives $\bigcup_{x \in Z} I_x \in \I(Q)$. 
\end{remark}

\begin{theorem}[Novak--Bandelt--Ern\'e]\label{thm:unioncomplete}
If $\Z$ is union-complete, then every $\Z$-continuous poset has the $\Z$-interpolation property, hence is strongly $\Z$-continuous. 
\end{theorem}

The original proof can be found in Novak \cite[Prop.~1.22]{Novak82b} and Bandelt and Ern\'e \cite[Theorem]{Bandelt83}. 
Yet, as these authors supposed $\Z$-completeness of the poset, we provide a proof without this hypothesis for the sake of rigor. 

\begin{proof}
Let $P$ be a $\Z$-continuous poset, and let $x, y \in P$ with $x \ll^{\Z} y$. 
By $\Z$-continuity, there is some $\Z$-subset $Z$ of $P$ included in $\twoheaddownarrow^{\Z} y$ with sup $y$. 
If $z \in P$, then by $\Z$-continuity $\twoheaddownarrow^{\Z} z$ is a $\Z$-ideal with sup equal to $z$. 
Now the subset ${ \twoheaddownarrow^{\Z} Z } := { \bigcup_{z \in Z} \twoheaddownarrow^{\Z} z }$ is a $\Z$-ideal of $P$ by union-completeness of $\Z$ and Remark~\ref{rk:unionc}. 
So let $Y$ be a $\Z$-subset of $P$ such that ${ \downarrow\!\! Y } = { \twoheaddownarrow^{\Z} Z }$. 
We show that the sup of $Y$ is $y$. 
It is obvious that $y$ is an upper bound of $Y$. 
Let $u \in P$ be another upper bound of $Y$. 
If $z \in Z$, then $\twoheaddownarrow^{\Z} z \subseteq { \twoheaddownarrow^{\Z} Z } = { \downarrow\!\! Y } \subseteq { \downarrow\!\! u }$. 
Thus, $z = \bigvee \twoheaddownarrow^{\Z} z \leqslant u$, for all $z \in Z$. 
This yields $y = \bigvee Z \leqslant u$. 
So we have shown that $y = \bigvee Y$. 
Now, $x \ll^{\Z} y = \bigvee Y$ and $Y$ is a $\Z$-subset of $P$, so $x \leqslant w$ for some $w \in Y$. 
Then $w \in \twoheaddownarrow^{\Z} Z$, so $w \ll^{\Z} v$ for some $v \in Z$. 
This proves that $x \ll^{\Z} v \ll^{\Z} y$. 
\end{proof}

For more on interrelations between union-completeness and the interpolation property, we refer the reader to Ern\'e \cite[Section~5]{Erne99}. 
For the special case $\Z = \D$ of classical domain theory, it is well known that every $\D$-continuous poset is $\D$-interpolating (so is strongly $\D$-continuous), see e.g.\ \cite[Theorem~I-1.9]{Gierz03}. 
However, the converse is false in general, as testified by the following examples. 

A \textit{complete lattice} is a poset in which every subset has a sup. Thus, $\A$-complete posets are the same objects as complete lattices. 

\begin{figure}
	\centering
	\begin{tikzpicture}%[scale=.7]
		\node (one) at (0,4*1) {$1$};
		\node (a) at (0,4*1/2) {$1/2$};
		\node (b) at (0, 4*3/4) {$3/4$};
		%\node (c) at (0, 4*7/8) {$7/8$};
		\node (zero) at (0,0) {$0$};
		\node (o) at (3,4*1/2) {$\omega$};
		\draw[thick] (zero) -- (a) -- (b);% -- (c);
		\draw[thick] (one) -- (o) -- (zero);
		\draw[thick, dotted] (b) -- (one);%(c) -- (one);
	\end{tikzpicture}
	\caption{Hasse diagram of Example~\ref{ex:one} }
	\label{fig:Hasse}
\end{figure}

\begin{example}\label{ex:one}
Figure~\ref{fig:Hasse} depicts the Hasse diagram of the countable complete lattice $P = \{ 0, 1/2, 3/4, \ldots, 1 \} \cup \{ \omega \}$, where $\omega$ is an element comparable only with $0$ and $1$, and such that $0 < \omega < 1$. 
The poset $P$ has the $\D$-interpolation property.  
Indeed, $0 \ll^{\D} \omega$ and $x \ll^{\D} y \ll^{\D} 1$ for all $x, y \in P \!\setminus\! \{ \omega, 1 \}$ such that $x \leqslant y$, and no other $\D$-below relation holds. 
However, $P$ is not $\D$-continuous since $\twoheaddownarrow^{\D} \omega = \{ 0 \}$; indeed, we have $0 \ll^{\D} \omega$ but not $\omega \ll^{\D} \omega$, because $\omega \leqslant \bigvee (P \!\setminus\! \{ \omega, 1 \}) = 1$. 
\end{example}

\begin{figure}
	\centering
	\begin{tikzpicture}%[scale=.7]
		\node (z) at (1,0) {$0$};
		\node (zu) at (0,1) {$1$};
		\node (zd) at (0,2) {$2$};
		\node (zt) at (0,3) {$3$};
		\node (uu) at (2,1) {$1'$};
		\node (ud) at (2,2) {$2'$};
		\node (ut) at (2,3) {$3'$};
		\node (o) at (1,4) {$\omega$};
		\draw[thick] (z) -- (zu) -- (zd) -- (zt);
		\draw[thick] (z) -- (uu) -- (ud) -- (ut);
		\draw[thick, dotted] (zt) -- (o) -- (ut);
	\end{tikzpicture}
	\caption{Hasse diagram of Example~\ref{ex:two} }
	\label{fig:Hasse2}
\end{figure}

\begin{example}\label{ex:two}
The poset defined by the Hasse diagram of Figure~\ref{fig:Hasse2} is the disjoint union of two copies of $\mathbb{N}\setminus\{0\}$, with a bottom element $0$ and a top element $\omega$ added. It is a countable complete lattice with the $\D$-interpolation property, but it is not $\D$-continuous since $\twoheaddownarrow^{\D} \omega = \{ 0 \}$. 
\end{example}

\begin{problem}
It is a well-known result of order theory that a lattice is not distributive if and only if it contains a diamond (the $5$-element lattice called $M_3$) or a pentagon (the $5$-element lattice called $N_5$). 
In the same spirit, is it possible to find a general pattern for $\D$-interpolating posets that are not $\D$-continuous? 
Or is it at least possible to exhibit a large class of such posets? 
\end{problem}

\subsection{Bases and algebraicity}

In the literature, a \textit{$\Z$-basis} of a poset $P$ is usually defined as a subset $B$ of $P$ such that $\twoheaddownarrow^{\Z} x \cap B$ contains a $\Z$-subset of $B$ whose sup in $P$ is $x$, for all $x \in P$. See e.g.\ Venugopalan \cite[Definition~3.1]{Venugopalan86}. 
In this paper, we shall use a somewhat weaker notion. 

\begin{definition}
A subset $B$ of a poset $P$ is a \textit{weak $\Z$-basis} of $P$ if $\downarrow\!\!(\twoheaddownarrow^{\Z} x \cap B)$ contains a $\Z$-subset of $P$ whose sup in $P$ is $x$, for all $x \in P$. 
\end{definition}

So, given $x \in P$, the $\Z$-subset of the definition of a weak $\Z$-basis $B$ is included in $\downarrow\!\!(\twoheaddownarrow^{\Z} x \cap B)$, not in $\twoheaddownarrow^{\Z} x \cap B$, and is in $\Z(P)$, not in $\Z(B)$. 

\begin{definition}\label{def:dense}
A subset $B$ of a poset $P$ is \textit{$\Z$-dense in $P$} (resp.\ \textit{strongly $\Z$-dense in $P$}) if, whenever $x \ll^{\Z} y$, there is some $b \in B$ with $x \leqslant b \ll^{\Z} y$ (resp.\ $x \ll^{\Z} b \ll^{\Z} y$). 
\end{definition}

Note that a poset is $\Z$-continuous if and only if it has a (weak) $\Z$-basis, and $\Z$-interpolating if and only if it has a strongly $\Z$-dense subset.  
Moreover, if $P$ is a $\Z$-interpolating poset with $\Z$-dense subset $B$, then $B$ is strongly $\Z$-dense. 

\begin{proposition}\label{prop:basis}
Let $P$ be a poset and $B \subseteq P$. 
Consider the following assertions: 
\begin{enumerate}
  \item\label{prop:basis1} $B$ is a $\Z$-basis of $P$;
  \item\label{prop:basis2} $B$ is a weak $\Z$-basis of $P$;
  \item\label{prop:basis3} $P$ is $\Z$-continuous and $B$ is $\Z$-dense in $P$. 
\end{enumerate}
Then \eqref{prop:basis1} $\Rightarrow$ \eqref{prop:basis2} $\Leftrightarrow$ \eqref{prop:basis3}.
Moreover, if $\Z \in \{ \A, \A^*, \D, \F \}$, then all three assertions are equivalent. 
\end{proposition}

\begin{proof}
\eqref{prop:basis1} $\Rightarrow$ \eqref{prop:basis2} is clear from the definitions. 

\eqref{prop:basis2} $\Rightarrow$ \eqref{prop:basis3}: 
Since $P$ has a weak $\Z$-basis, it is $\Z$-continuous. 
Let $x \ll^{\Z} y$. 
By \eqref{prop:basis2} there is a $\Z$-subset $Z$ of $P$ included in $\downarrow\!\!(\twoheaddownarrow^{\Z} y \cap B)$ such that $y = \bigvee Z$. 
Then $x \in { \downarrow\!\! Z }$, so $x \in { \downarrow\!\!(\twoheaddownarrow^{\Z} y \cap B) }$. 
Thus, $x \leqslant b$ for some $b \in { \twoheaddownarrow^{\Z} y \cap B }$. 
So $x \leqslant b \ll^{\Z} y$, as required.

\eqref{prop:basis3} $\Rightarrow$ \eqref{prop:basis2}: 
Let $y \in P$. 
By $\Z$-continuity of $P$, there is some $\Z$-subset $Z$ of $P$ included in $\twoheaddownarrow^{\Z} y$ with sup $y$. 
Now $B$ is $\Z$-dense in $P$, so ${ \twoheaddownarrow^{\Z} y } \subseteq { \downarrow\!\!(\twoheaddownarrow^{\Z} y \cap B) }$, so that ${ \twoheaddownarrow^{\Z} y } = { \downarrow\!\!(\twoheaddownarrow^{\Z} y \cap B) }$. 
This proves that $B$ is a weak $\Z$-basis.

\eqref{prop:basis3} $\Rightarrow$ \eqref{prop:basis1} if $\Z \in \{ \A, \A^*, \D, \F \}$:  
Let $x \in P$. 
Since $P$ is $\Z$-continuous, there is some $\Z$-subset $Z$ of $P$ with sup $x$ such that ${ \twoheaddownarrow^{\Z} x } = { \downarrow\!\! Z }$. 
If $z \in Z$, then by \eqref{prop:basis3} there is some $b_z \in B$ such that $z \leqslant b_z \ll^{\Z} x$. 
Take $Z' := { \{ b_z : z \in Z \} }$. 
Then $Z \subseteq { \downarrow\!\! Z' }$ on the one hand, 
and $Z' \subseteq { \twoheaddownarrow^{\Z} x \cap B } \subseteq { \twoheaddownarrow^{\Z} x } = { \downarrow\!\! Z }$ on the other hand, so that ${ \downarrow\!\! Z' } = { \downarrow\!\! Z }$. 
This implies that $Z'$ has a sup equal to $x$. 
Now, if $\Z \in \{ \A, \A^*, \F \}$, then $Z'$ is obviously a $\Z$-subset of $B$; and if $\Z = \D$, then ${ \downarrow\!\! Z' } = { \downarrow\!\! Z }$ and $Z$ directed imply that $Z'$ is a directed subset of $B$. 
So $Z'$ is a $\Z$-subset of $B$ included in $\twoheaddownarrow^{\Z} x \cap B$ with sup $x$. 
Thus, $B$ is a $\Z$-basis of $P$. 
%
%Since $P$ is $\D$-continuous, there is some directed subset $D$ of $P$ with sup $x$ such that ${ \twoheaddownarrow^{\D} x } = { \downarrow\!\! D }$. 
%If $d \in D$, then by \eqref{prop:basis3} there is some $b_d \in B$ such that $d \leqslant b_d \ll^{\D} x$. 
%Take $D' := { \{ b_d : d \in D \} }$. 
%Then $D \subseteq { \downarrow\!\! D' }$ on the one hand, and $D' \subseteq { \twoheaddownarrow^{\D} x \cap B } \subseteq { \twoheaddownarrow^{\D} x } = { \downarrow\!\! D }$ on the other hand, so that ${ \downarrow\!\! D' } = { \downarrow\!\! D }$. 
%This implies that $D'$ is directed and has a sup equal to $x$. 
%So $D'$ is a directed subset of $B$ included in $\twoheaddownarrow^{\D} x \cap B$ with sup $x$. 
%This proves that $B$ is a $\D$-basis of $P$. 
%The cases where $\Z \in \{ \A, \A^*, \F \}$ can be proved along similar (yet simpler) lines.
\end{proof}

A poset is \textit{$\Z$-inductive} (resp.\ \textit{$\Z$-algebraic}) if the subset of its $\Z$-compact elements forms a $\Z$-basis (resp.\ a weak $\Z$-basis). 
From the previous proposition we deduce that every $\Z$-inductive poset is $\Z$-algebraic and that every $\Z$-algebraic poset is strongly $\Z$-continuous. 

These definitions agree with Wright et al.\ \cite{Wright78} and Bandelt and Ern\'e \cite{Bandelt83}, except that the latter authors included $\Z$-completeness in the definition of $\Z$-algebraicity. 
Yet, the example of the complete chain $\omega + 1 = \{ 0, 1, \ldots, \omega \}$ that they provided remains a valid instance of a $\B$-algebraic poset that is not $\B$-inductive: as they explained, there is no subset made of $\B$-compact elements that has a $\B$-compact upper bound and admits $\omega$ as its sup. 
It shows the non-equivalence of the notions of $\Z$-basis and weak $\Z$-basis in general. 

%As shown by the following result, this definition is equivalent to the usual one \cite[Definition~III-4.1]{Gierz03}, yet it will prove more convenient in the frame of this paper. 

%\begin{proposition}
%A subset $B$ of a poset $P$ is a $\Z$-basis of $P$ if and only if the subset $\twoheaddownarrow^{\Z} x \cap B$ contains a $\Z$-subset whose sup is $x$ for all $x \in P$. 
%\end{proposition}

%\begin{proof}
%The `if' part is straightforward, we prove the `only if' part. 
%Let $B$ be a $\Z$-basis of $P$ and let $x \in P$. 
%By hypothesis there is some $\Z$-subset $Z$ of $\downarrow\!\! (\twoheaddownarrow^{\Z} x \cap B)$ whose sup is $x$. 
%Consequently, $\twoheaddownarrow^{\Z} x \cap B$ is nonempty and also has a sup equal to $x$. 
%We prove that $\twoheaddownarrow^{\Z} x \cap B$ is directed. 
%Let $b_1, b_2 \in \twoheaddownarrow x \cap B$. 
%The poset $P$ is continuous, so $\twoheaddownarrow x$ is directed, thus there is some $y \ll x$ such that 
%$b_1 \leqslant y$ and $b_2 \leqslant y$. 
%Now $x = \bigvee D$, hence $y \in \downarrow\!\! D \subseteq \downarrow\!\! (\twoheaddownarrow x \cap B)$. 
%This entails the existence of some $b \in B$ such that $y \leqslant b \ll x$. 
%Thus, $b_1 \leqslant b$ and $b_2 \leqslant b$, which shows that $\twoheaddownarrow x \cap B$ is directed. 
%\end{proof}

\subsection{$\Z$-sup-preserving maps and Galois connections}\label{subsec:sl}

A map $f : P \to Q$ between posets is \textit{$\Z$-sup-preserving} if, for every $\Z$-subset $Z$ of $P$ with sup in $P$, the sup of $f(Z)$ exists in $Q$ and $f(\bigvee_P Z) = \bigvee_Q f(Z)$.  
In the classical case, a $\D$-sup-preserving map is the same as a \textit{Scott-continuous} map. 
This is why, in the literature, a $\Z$-sup-preserving map has been actually called a \textit{$\Z$-continuous} map; however, for the sake of clarity, we shall dedicate the term \textit{$\Z$-continuous} to posets in this paper. 
Recall from \cite[Corollary~1.7]{Baranga96} that every $\Z$-sup-preserving map is order-preserving, thanks to Condition $iii)$ in the definition of a proper subset system. 
Note also that $\Z$-sup-preserving maps are stable under composition. 

The following lemma will be used repeatedly along this paper. 

\begin{lemma}\label{lem:image}
If a map $f : P \to Q$ is $\Z$-sup-preserving, then its corestriction $f^{\circ} : P \to f(P)$ is $\Z$-sup-preserving. 
\end{lemma}

\begin{proof}
Let $Z$ be a $\Z$-subset of $P$ with sup. 
Since $f$ is $\Z$-sup-preserving, $f(Z)$ has a sup in $Q$ and $\bigvee_Q f(Z) = f(\bigvee_P Z) \in f(P)$. 
Thus, $f(\bigvee_P Z)$ is also the sup of $f(Z)$ in $f(P)$, and $f^{\circ}(\bigvee_P Z) = \bigvee_{f(P)} f^{\circ}(Z)$. 
\end{proof}

A \textit{Galois connection} between two posets $P$ and $Q$ is a pair $(f, g)$ of order-preserving maps $f : P \to Q$ and $g : Q \to P$ such that $f(x) \leqslant y$ if and only if $x \leqslant g(y)$, for all $x \in P$ and $y \in Q$. 
The map $f$ (resp.\ $g$) is the \textit{left adjoint} (resp.\ \textit{right adjoint}) of the Galois connection. 
We refer the reader to Ern\'{e} et al.\ \cite{Erne93} for usual properties and various references on Galois connections. 
Note the following properties: 
\begin{itemize}
  \item $x \leqslant g(f(x))$, for all $x \in P$;
  \item $f(g(y)) \leqslant y$,	 for all $y \in Q$; 
  \item $g(f(x)) = x$ for all $x \in P$ iff $f$ is injective iff $g$ is surjective;
  \item $f(g(y)) = y$ for all $y \in Q$ iff $g$ is injective iff $f$ is surjective;
  \item $f \circ g \circ f = f$ and $g \circ f \circ g = g$;
  \item $f(P) = f(g(Q))$ and $g(Q) = g(f(P))$. 
\end{itemize} 

We also recall the following well-known result. 

\begin{lemma}\label{lem:galois}
The left adjoint of a Galois connection is $\A$-sup-preserving (hence $\Z$-sup-preserving for every proper subset system $\Z$). 
\end{lemma}

\begin{proof}
See e.g.\ \cite[Proposition~7(3)]{Erne93}. 
\end{proof}

We will often generalize results involving Galois connections to pre-Galois connections, that we define as follows. 
A \textit{pre-Galois connection} is a pair $(f, g)$ of order-preserving maps $f : P \to Q$ and $g : Q \to P$ such that $f \circ g \circ f = f$ and $f(g(y)) \leqslant y$, for all $y \in Q$. 
The map $f$ (resp.\ $g$) is the \textit{left adjoint} (resp.\ \textit{right adjoint}) of the pre-Galois connection.

%%%%%%%%%%%%%%%%%%%%%%
%%%%%%%%%%%%%%%%%%%%%%
%%%%%%%%%%%%%%%%%%%%%%
%%%%%%%%%%%%%%%%%%%%%%
\section{The $\Z$-refinement property and $\Z$-adequacy}\label{sec:maps}

In this section, we focus on order-preserving maps between posets and introduce the property of $\Z$-refinement.  
This property suitably replaces preservation of the $\Z$-below relation in cases where $\Z$-continuity is not given. 
%We show how this property is linked with other properties, in the line of Bandelt and Ern\'e \cite{Bandelt83}; 
%we also take the opportunity to recall and synthesize various results due to Venugopalan. 

An \textit{order-embedding} is an order-preserving map $f : P \to P'$ such that $x \leqslant_{P} y$ whenever $f(x) \leqslant_{P'} f(y)$, for all $x, y \in P$. 
Equivalently, the corestriction $f^{\circ}$ is an isomorphism of posets. %injective and is the right adjoint of a Galois connection. 
In particular, an order-embedding is always injective and quasi-invertible. 
Note also that every injective left adjoint of a Galois connection is an order-embedding. 
An order-preserving map $f : P \to P'$ \textit{reflects the $\Z$-below relation} if $f(x) \ll^{\Z}_{P'} f(y)$ implies $x \ll^{\Z}_P y$, for all $x, y \in P$.

The following lemma generalizes \cite[Lemma~1.2]{Hoffmann85}. 

\begin{lemma}\label{lem:reflecting}
Let $f : P \to P'$ be a $\Z$-sup-preserving order-embedding. 
Then $f$ reflects the $\Z$-below relation. 
\end{lemma}

\begin{proof}
Suppose that $f(x) \ll^{\Z}_{P'} f(y)$ for some $x, y \in P$. 
Let $Z$ be a $\Z$-subset of $P$ with sup such that $y \leqslant \bigvee_P Z$. 
The map $f$ preserves $\Z$-sups, so $f(y) \leqslant f(\bigvee_P Z) = \bigvee_{P'} f(Z)$. 
Since $f(Z)$ is a $\Z$-subset of $P'$ and $f(x) \ll^{\Z}_{P'} f(y)$, there exists some $z \in Z$ such that $f(x) \leqslant f(z)$. 
Since $f$ is an order-embedding by hypothesis, this yields $x \leqslant z$. 
This proves that $x \ll^{\Z}_P y$. 
\end{proof}

The following result will not be used in this paper, but is interesting in its own right. 

\begin{proposition}
Let $P'$ be a $\Z$-continuous poset, and let $f : P \to P'$ be a surjective order-preserving map that reflects the $\Z$-below relation. 
Then $f$ is $\Z$-sup-preserving. 
\end{proposition}

\begin{proof}
Let $Z$ be a $\Z$-subset of $P$ with sup, and take $x' := f(\bigvee_P Z)$. 
Since $P'$ is $\Z$-continuous, there exists some $\Z$-subset $Z'$ of $P'$ included in $\twoheaddownarrow_{P'} x'$ such that $x' = \bigvee_{P'} Z'$. 
Using the properties on $f$, it is not difficult to show that ${ Z' } \subseteq { \downarrow_{P'}\!\! f(Z) } \subseteq { \downarrow_{P'}\!\! x' }$. 
Consequently, $f(Z)$ also admits $x'$ as its sup in $P'$, i.e.\ $f(\bigvee_P Z) = x' = \bigvee_{P'} f(Z)$. 
This proves that $f$ is $\Z$-sup-preserving. 
\end{proof}

\begin{definition}
An order-preserving map $f : P \to P'$ has the \textit{$\Z$-refinement property} if, whenever $x \in P$ and $f(x) \leqslant \bigvee_{P'} Z'$ for some $\Z$-subset $Z'$ of $P'$ with sup, there exists a $\Z$-subset $Z$ of $P$ with sup such that $f(Z) \subseteq { \downarrow_{P'}\!\! Z' }$ and $x \leqslant \bigvee_P Z$. 
An order-preserving map is \textit{$\Z$-adequate} if it is $\Z$-sup-preserving and has the $\Z$-refinement property. 
\end{definition}

Note that order-preserving maps with the $\Z$-refinement property (resp.\ $\Z$-adequate maps) are stable under composition. 

We say that a map $f : P \to P'$ \textit{preserves the $\Z$-below relation} if $x \ll^{\Z}_P y$ implies $f(x) \ll^{\Z}_{P'} f(y)$, for all $x, y \in P$. 
The following result shows the close relation between the $\Z$-refinement property and the preservation of the $\Z$-below relation. 

\begin{proposition}\label{prop:refinement}
Let $f : P \to P'$ be an order-preserving map. 
If $f$ has the $\Z$-refinement property, then $f$ preserves the $\Z$-below relation. 
Conversely, if $P$ is $\Z$-continuous and $f$ preserves the $\Z$-below relation, then $f$ has the $\Z$-refinement property. 
\end{proposition}

\begin{proof}
First, suppose that $f$ has the $\Z$-refinement property, and let $x \ll^{\Z}_P y$ and $Z'$ be a $\Z$-subset of $P'$ with sup such that $f(y) \leqslant \bigvee_{P'} Z'$. 
By the $\Z$-refinement property there exists a $\Z$-subset $Z$ of $P$ with sup such that $f(Z) \subseteq { \downarrow_{P'}\!\! Z' }$ and $y \leqslant \bigvee_P Z$. 
Since $x \ll^{\Z}_P y$ we have $x \in { \downarrow_{P}\!\! Z }$, so $f(x) \in f(\downarrow_{P}\!\! Z) \subseteq { \downarrow_{P'}\!\! f(Z) } \subseteq { \downarrow_{P'}\!\! Z' }$. 
This proves that $f(x) \ll^{\Z}_{P'} f(y)$. 

Conversely, suppose that $P$ is $\Z$-continuous and $f$ preserves the $\Z$-below relation. 
Let $f(x) \leqslant \bigvee_{P'} Z'$ where $x \in P$ and $Z'$ is a $\Z$-subset of $P'$ with sup. 
Then $\twoheaddownarrow^{\Z}_{P'} f(x) \subseteq { \downarrow_{P'}\!\! Z' }$ by definition of the $\Z$-below relation. 
Since $f$ preserves the $\Z$-below relation, we have $f(\twoheaddownarrow^{\Z}_P x) \subseteq { \twoheaddownarrow^{\Z}_{P'} f(x) }$, so $f(\twoheaddownarrow^{\Z}_P x) \subseteq { \downarrow_{P'}\!\! Z' }$. 
Moreover, $P$ is $\Z$-continuous, so there is some $\Z$-subset $Z$ of $P$ included in $\twoheaddownarrow^{\Z}_P x$ whose sup is $x$. 
Thus, $f(Z) \subseteq { \downarrow_{P'}\!\! Z' }$ and $x = \bigvee_P Z$. 
This shows that $f$ has the $\Z$-refinement property. 
\end{proof}

The following extends a result of Bandelt and Ern\'e \cite{Bandelt83}. 

\begin{proposition}\label{prop:left}
The left adjoint of a Galois connection has the $\Z$-refin\-ement property if and only if its right adjoint is $\Z$-sup-preserving. 
In this case, it is $\Z$-adequate and preserves the $\Z$-below relation. 
\end{proposition}

\begin{proof}
Let $f : P \to P'$ be the left adjoint and $g : P' \to P$ be the right adjoint of a Galois connection. 
First assume that $f$ has the $\Z$-refinement property. 
Let $Z'$ be a $\Z$-subset of $P'$ with sup, and let $u \in P$ be an upper bound of $g(Z')$. 
We take $x := g(\bigvee_{P'} Z')$. 
Then $f(x) \leqslant \bigvee_{P'} Z'$, so by the $\Z$-refinement property there exists a $\Z$-subset $Z$ of $P$ with sup such that $f(Z) \subseteq { \downarrow_{P'}\!\! Z' }$ and $x \leqslant \bigvee_P Z$. 
The former property for $Z$ is equivalent to $Z \subseteq { \downarrow_{P}\!\! g(Z') }$, 
so $g(\bigvee_{P'} Z') = x \leqslant \bigvee_P Z \leqslant u$. 
This shows that $g(\bigvee_{P'} Z')$ is the sup of $g(Z')$, so $g$ is $\Z$-sup-preserving. 

Conversely, assume that $g$ is $\Z$-sup-preserving. 
Let $Z'$ be a $\Z$-subset of $P'$ with sup such that $f(x) \leqslant \bigvee_{P'} Z'$. 
Then $x \leqslant g(\bigvee_{P'} Z') = \bigvee_P g(Z')$, since $g$ is $\Z$-sup-preserving. 
Moreover, we have $f(g(Z')) \subseteq { \downarrow_{P'}\!\! Z' }$. 
This proves that $f$ has the $\Z$-refinement property.  
\end{proof}

\begin{definition}
A subset $A$ of a poset $P$ \textit{$\Z$-sup-generates} $P$, or is \textit{$\Z$-sup-generating} in $P$ if, for every $x \in P$, there is a $\Z$-subset $Z$ of $P$ included in $A$ with sup $x$. 
\end{definition}

Note that a $\Z$-sup-generating subset must contain all $\Z$-compact elements. 

\begin{example}\label{ex:supdense}
In a poset $P$, let $C$ be the subset of $P$ made of its $\Z$-compact elements. 
Consider the following assertions.
\begin{enumerate}
  \item\label{ex:supdense1} $P$ is $\Z$-inductive; 
  \item\label{ex:supdense2} $C$ $\Z$-sup-generates $P$;
  \item\label{ex:supdense3} $P$ is $\Z$-algebraic. 
\end{enumerate}
Then \eqref{ex:supdense1} $\Rightarrow$ \eqref{ex:supdense2} $\Rightarrow$ \eqref{ex:supdense3}.  
Indeed, it is easily seen that $C$ $\Z$-sup-generates $P$ if and only if, for every $x \in P$, there is some $\Z$-subset $Z$ of $P$ included in $\twoheaddownarrow^{\Z} x \cap C$ whose sup in $P$ is $x$. 
The rest follows from the definitions.
\end{example}

\begin{theorem}\label{thm:imageofzcp}
Let $P$ be a $\Z$-continuous poset with weak $\Z$-basis $B$, and let $f : P \to P'$ be a $\Z$-adequate  quasi-invertible map such that $f(P)$ $\Z$-sup-generates $P'$. 
Moreover, assume that $\Z$ is union-complete. %one of the following conditions holds: 
%\begin{enumerate}
%  \item\label{thm:imageofzcp1} $f(P)$ is $\Z$-sup-generating in $P'$, $f$ is quasi-invertible, and $\Z$ is union-complete;
%  \item\label{thm:imageofzcp2} $f$ is surjective. 
%\end{enumerate}
Then $P'$ is strongly $\Z$-continuous with weak $\Z$-basis $f(B)$. 
\end{theorem}

\begin{proof}
%Case \eqref{thm:imageofzcp1}:
Let $g : f(P) \to P$ be a quasi-inverse of $f$, and let $x' \in P'$. 
Since $f(P)$ is $\Z$-sup-generating in $P'$, there is a $\Z$-subset $Z'$ of $P'$ included in $f(P)$ such that $x' = \bigvee_{P'} Z'$. 
By Proposition~\ref{prop:stable}\eqref{prop:stable7}, $Z'$ is a $\Z$-subset of $f(P)$. 
For each $y' \in f(P)$, take $x_{y'} := g(y') \in P$. 
We have $y' = f(x_{y'})$, and ${ I_{y'} } := { \twoheaddownarrow^{\Z}_P x_{y'} } = { \downarrow_{P}\!\! (\twoheaddownarrow^{\Z}_P x_{y'} \cap B) }$ is a $\Z$-ideal of $P$ with sup $x_{y'}$. 
This implies that $\downarrow_{P'}\!\! f(I_{y'})$ is a $\Z$-ideal of $P'$. 
Consider the subset ${ I' } := { \bigcup_{z' \in Z'} \downarrow_{P'}\!\! f(I_{z'}) }$. 
Using union-completeness of $\Z$, we can apply Theorem~\ref{thm:unioncomplete} together with Remark~\ref{rk:unionc}, which yields $I' \in \I(P')$. 
Moreover, 
\begin{align*}
{ I' } &= { \bigcup_{z' \in Z'} \downarrow_{P'}\!\! f(\downarrow_{P}\!\! (\twoheaddownarrow^{\Z}_P x_{z'} \cap B)) } \\
&\subseteq { \bigcup_{z' \in Z'} \downarrow_{P'}\!\! f(\twoheaddownarrow^{\Z}_P x_{z'} \cap B) } = { \downarrow_{P'}\!\! \bigcup_{z' \in Z'} f(\twoheaddownarrow^{\Z}_P x_{z'} \cap B) } \\
&\subseteq { \downarrow_{P'}\!\! \bigcup_{z' \in Z'} \left( f(\twoheaddownarrow^{\Z}_P x_{z'}) \cap f(B) \right) }. 
\end{align*}
Since $f$ has the $\Z$-refinement property, it preserves the $\Z$-below relation by Proposition~\ref{prop:refinement}, which implies that $f(\twoheaddownarrow^{\Z}_P x_{z'}) \subseteq { \twoheaddownarrow^{\Z}_{P'} f(x_{z'}) } = { \twoheaddownarrow^{\Z}_{P'} z' }$, for all $z' \in Z'$. 
Hence, 
\[
I' \subseteq { \downarrow_{P'}\!\! \bigcup_{z' \in Z'} \left( \twoheaddownarrow^{\Z}_{P'} z' \cap f(B) \right) } \subseteq { \downarrow_{P'}\!\! (\twoheaddownarrow^{\Z}_{P'} x' \cap f(B)) }. 
\]
In addition, the subset $I'$ has sup $x'$ in $P'$ (use the fact that $f$ is $\Z$-sup-preserving). 
So $P'$ is $\Z$-continuous with weak $\Z$-basis $f(B)$. 
To conclude the proof, recall that $\Z$ is assumed to be union-complete, so $P'$ is strongly $\Z$-continuous by Theorem~\ref{thm:unioncomplete}. 
%Case~\eqref{thm:imageofzcp1} is straightforward, for in this case $\Z$ is supposed to be union-complete: the $\Z$-interpolation property automatically follows from Theorem~\ref{thm:unioncomplete}. 
%Let us consider Case~\eqref{thm:imageofzcp2}. 
%Let $x', y' \in P'$ such that $x' \ll^{\Z}_{P'} y'$. 
%Since $f$ is surjective, there is some $y \in P$ with $y' = f(y)$. 
%Moreover, $P$ is $\Z$-continuous, so there is some $\Z$-subset $Z$ of $P$ included in $\twoheaddownarrow_P^{\Z} y$ with sup equal to $y$. 
%Now, $f$ is $\Z$-sup-preserving, so 
%\[
%x' \ll^{\Z}_{P'} y' = f(y) = f(\bigvee_P Z) = \bigvee_{P'} f(Z).
%\] 
%Thus, there is some $z \in Z \subseteq { \twoheaddownarrow_P^{\Z} y }$ such that $x' \leqslant f(z)$. 
%We have supposed that $P$ is $\Z$-interpolating, so $z \ll_P^{\Z} t \ll_P^{\Z} y$ for some $t \in P$. 
%Since $f$ has the $\Z$-refinement property, it preserves the $\Z$-below relation by Proposition~\ref{prop:refinement}, so $x' \leqslant f(z) \ll^{\Z}_{P'} f(t) \ll^{\Z}_{P'} f(y) = y'$. 
%This proves that $P'$ is $\Z$-interpolating. 
\end{proof}

%%%%%%%%%%%%%%%%%%%%%%
%%%%%%%%%%%%%%%%%%%%%%
%%%%%%%%%%%%%%%%%%%%%%
%%%%%%%%%%%%%%%%%%%%%%
\section{Sufficient conditions to be a $\Z$-subposet}\label{sec:main}

In this section, we specialize the definitions of the previous section to subsets of posets seen as inclusion maps, 
and we state and prove our Main Theorem. 
In passing, we also examine the question of transferring the continuity properties of a covering family of subposets of a poset to the poset itself. 

\subsection{$\Z$-subposets and their properties}

Let $P$ be a poset. 
A subset $A$ of $P$ is called \textit{$\Z$-sup-preserving} if every $\Z$-subset $Z$ of $A$ with sup in $A$ also has a sup in $P$, and $\bigvee_P Z = \bigvee_A Z$; this amounts to saying that the inclusion map $A \to P$ is $\Z$-sup-preserving. 

\begin{example}
\begin{itemize}
  \item Every chain in a poset is $\F$-sup-preserving. 
  \item Every finite subset of a poset is $\D$-sup-preserving. 
  \item Every antichain in a poset is $\A$-sup-preserving, where an \textit{antichain} is a subset $A$ such that $x \not\leqslant y$ for all $x, y \in A$ with $x \neq y$. 
  \item Every poset is $\A$-sup-preserving in its Dedekind--MacNeille completion (see, e.g.\ \cite{Erne93} for a definition). 
\end{itemize}
\end{example}

%The following example is given by Baranga \cite[Lemma~1.11]{Baranga96}. 

%\begin{example}[Baranga]\label{ex:baranga}
%Let $P$ be a conditionally $\Z$-complete poset. 
%If $p : P \to P$ is a $\Z$-sup-preserving projection, then $p(P)$ is a $\Z$-sup-preserving subset of $P$. 
%\end{example}

A poset is \textit{conditionally $\Z$-complete} if every upper-bounded $\Z$-subset has a sup. 
We sometimes use the abbreviation \textit{cond.\ $\Z$-complete}.  
Conditionally $\A^*$-complete posets are thus the same object as conditionally complete posets; and conditionally $\A$-complete nonempty posets as conditionally complete posets with a least element. 

Case~\eqref{prop:baranga4} of the following proposition was given by \cite[Lemma~1.11]{Baranga96}. 

\begin{proposition}\label{prop:baranga}
Let $f : P \to P'$ be an order-preserving map. 
Assume that one of the following conditions holds: 
\begin{enumerate}
  \item\label{prop:baranga0} $f$ is surjective; 
  \item\label{prop:baranga1} $f$ is the left adjoint of a Galois connection; 
  \item\label{prop:baranga1bis} $f$ is the left adjoint of a pre-Galois connection; 
  \item\label{prop:baranga4} $P$ is cond.\ $\Z$-complete and $f$ is a $\Z$-sup-preserving projection;
  \item\label{prop:baranga2} $P$ is cond.\ $\Z$-complete and $f$ is $\Z$-sup-preserving quasi-invertible; 
  \item\label{prop:baranga2bis} $f$ is $\Z$-sup-preserving and has a $\Z$-sup-preserving quasi-inverse; 
  \item\label{prop:baranga0bis} $f$ is a $\Z$-sup-preserving order-embedding;
  \item\label{prop:baranga7} $P'$ is cond.\ $\Z$-complete, $f$ is $\Z$-sup-preserving, and there is a $\Z$-sup-preserving map $g : P' \to P$ such that $f \circ g \circ f = f$. 
\end{enumerate}
Then $f(P)$ is $\Z$-sup-preserving in $P'$. 
\end{proposition}

\begin{proof}
Case \eqref{prop:baranga0} is trivial. 
Case \eqref{prop:baranga1} is a consequence of Case \eqref{prop:baranga1bis}. 
Case \eqref{prop:baranga4} is a consequence of Case \eqref{prop:baranga2}, since every projection is quasi-invertible. 

Case \eqref{prop:baranga2bis}: Let $g$ be a $\Z$-sup-preserving quasi-inverse of $f$. 
The injection map $i : f(P) \to P'$ can be written as $i = f \circ g$. 
Since both $f$ and $g$ are $\Z$-sup-preserving, so is $i$. 
Thus, $f(P)$ is $\Z$-sup-preserving in $P'$. 

Case \eqref{prop:baranga0bis}: Since $f$ is an order-embedding, $f^{\circ}$ is an isomorphism of posets with inverse map $j : f(P) \to P$, and the injection map $i : f(P) \to P'$ can be written as $i = f \circ j$. Now, $f$ is $\Z$-sup-preserving and $j$ is $\A$-sup-preserving (e.g.\ by Lemma~\ref{lem:galois}), so $i$ is $\Z$-sup-preserving, which means that $f(P)$ is $\Z$-sup-preserving. 

For Cases \eqref{prop:baranga1bis} and \eqref{prop:baranga2}, let $Z'$ be a $\Z$-subset of $f(P)$ with a sup in $f(P)$ denoted by $f(x_0)$. 
Let $u'$ be an upper bound of $Z'$ in $P'$. 
We take $Z := g(Z')$, where $g$ is either a right adjoint of $f$ in Case \eqref{prop:baranga1bis}, or a quasi-inverse of $f$ in Case \eqref{prop:baranga2}. 
Then $Z$ is a $\Z$-subset of $P$. 
Using $Z' \subseteq { f(P) }$ and the relation between $f$ and $g$, we get $f(Z) = Z'$. 

Case \eqref{prop:baranga1bis}:
We have $z' = f(g(z')) \leqslant f(g(u'))$, for all $z' \in Z'$. 
So $f(g(u'))$ is an upper bound of $Z'$ in $f(P)$, hence $f(x_0) \leqslant f(g(u'))$ by definition of $f(x_0)$. 
This yields $f(x_0) \leqslant u'$, so $f(x_0)$ is the sup of $Z'$ in $P'$. 

Case \eqref{prop:baranga2}:
Since $P$ is conditionally $\Z$-complete, $Z$ admits a sup $z_0$ in $P$. 
Now $f$ is supposed to be $\Z$-sup-preserving, so $f(z_0) = \bigvee_{P'} f(Z) = \bigvee_{P'} Z'$. 
So the sup of $Z'$ in $P'$ exists and belongs to $f(P)$, which implies that it equals $f(x_0)$. 

Case \eqref{prop:baranga7}: 
Since $P'$ is conditionally $\Z$-complete, $Z'$ admits a sup $z_0'$ in $P'$. 
Now $g$ is $\Z$-sup-preserving, so $Z$ has a sup in $P$ (equal to $g(z_0')$). 
For the rest of the proof, one can then follow the lines of the proof of Case \eqref{prop:baranga2}. 
%The map $f$ is also $\Z$-sup-preserving, so by Lemma~\ref{lem:image} $z_0' = \bigvee_{P'} Z' = \bigvee_{P'} f(Z) = f(\bigvee_{P} Z) = \bigvee_{f(P)} f(Z) = \bigvee_{f(P)} Z' = f(x_0)$. 
\end{proof}

\begin{definition}\label{def:refprop}
A subset $A$ of a poset $P$ has the \textit{$\Z$-refinement property} if the inclusion map $A \to P$ has the $\Z$-refinement property; equivalently, whenever $a \in A$ and $a \leqslant \bigvee_P Z$ for some $\Z$-subset $Z$ of $P$ with sup, there is a $\Z$-subset $Z'$ of $A$ included in $\downarrow_{P}\!\! Z$, with sup in $A$ such that $a \leqslant \bigvee_A Z'$. 
\end{definition}

\begin{example}
If $A$ is the set of $\Z$-compact elements of a poset $P$, then $A$ has the $\Z$-refinement property. 
Indeed, if $a$ is $\Z$-compact and $a \leqslant \bigvee_P Z$ for some $\Z$-subset $Z$, then $a \in { \downarrow_{P}\!\! Z }$. 
Taking $Z' = \{ a \}$, $Z'$ is a $\Z$-subset of $A$ included in $\downarrow_{P}\!\! Z$, and $a = \bigvee_A Z'$. 
\end{example}

\begin{definition}
A subset $A$ of a poset $P$ is a \textit{$\Z$-below preserving subposet} of $P$, or a \textit{$\Z$-subposet} of $P$ for short, if $A$ is $\Z$-sup-preserving in $P$ and the inclusion map $i : A \mapsto P$ preserves and reflects the $\Z$-below relation, that is, $x \ll^{\Z}_A y$ if and only if $x \ll^{\Z}_P y$, for all $x, y \in A$. 
\end{definition}

\begin{example}
In a poset, a singleton $\{ x \}$ is a $\Z$-subposet if and only if $x$ is a $\Z$-compact element. 
\end{example}

\begin{definition}\label{def:adequate}
A subset $A$ of a poset $P$ is \textit{$\Z$-adequate} if the inclusion map $i : A \to P$ is $\Z$-adequate, i.e.\ if $A$ is $\Z$-sup-preserving and has the $\Z$-refinement property in $P$. 
\end{definition}

\begin{example}\label{ex:adequate}
In a poset $P$, a subset $A$ that is both a lower set and an upper set is $\Z^*$-adequate. 
Indeed, consider first a nonempty $\Z$-subset $Z_1$ of $A$ with sup $a_1$ in $A$. 
If $u$ is an upper bound of $Z_1$ in $P$, then $u \in { \uparrow_{P}\!\! A }$ since $Z_1$ is nonempty, hence $u \in { A }$ since $A$ is an upper set. 
So $a_1 \leqslant u$ by definition of $a_1$. 
This yields $a_1 = \bigvee_P Z_1$. 
Hence, $A$ is $\Z^*$-sup-preserving. 
Now consider $a \in A$ and a $\Z$-subset $Z_2$ of $P$ with sup $x$ in $P$ such that $a \leqslant x$. 
Then $x \in { \uparrow_{P}\!\! A } = { A }$ and $Z_2 \subseteq { \downarrow_{P}\!\! x } \subseteq { \downarrow_{P}\!\! A} = { A }$. 
This implies that $\bigvee_A Z_2$ exists and equals $x = \bigvee_P Z_2$. 
Moreover, $Z_2$ is a $\Z$-subset of $A$ by Proposition~\ref{prop:stable}\eqref{prop:stable6}. 
This proves that $A$ has the $\Z$-refinement property. 
\end{example}

The following result extends Hoffmann \cite[Lemma~1.5]{Hoffmann85}. 
It notably gives a sufficient condition for a $\Z$-subposet to have the $\Z$-refinement property. 

\begin{theorem}\label{thm:main0}
Let $P$ be a poset and $A \subseteq P$. 
If $A$ is $\Z$-adequate, then $A$ is a $\Z$-subposet of $P$. 
Conversely, if $A$ is a $\Z$-continuous $\Z$-subposet of $P$, then $A$ is $\Z$-adequate. 
In the latter case, if moreover $\Z$ is union-complete and $A$ $\Z$-sup-generates $P$, then $P$ is $\Z$-continuous with weak $\Z$-basis $B$, where $B$ is any weak $\Z$-basis of $A$.  
\end{theorem}

%begin{corollary}
%Let $P$ be a poset and $A$ be a $\Z$-continuous $\Z$-subposet of $P$ that $\Z$-sup-generates $P$. 
%Moreover, assume that $\Z$ is union-complete. 
%Then $P$ is $\Z$-continuous with weak $\Z$-basis $B$, where $B$ is any weak $\Z$-basis of $A$.  
%\end{corollary}

%\begin{proof}
%Note that the inclusion map $i : A \to P$ is $\Z$-adequate by Theorem~\ref{thm:main0} and quasi-invertible with quasi-inverse the identity map of $A$. 
%Then apply Theorem~\ref{thm:imageofzcp} to $i$. 
%\end{proof}

\begin{proof}
Suppose first that $A$ is $\Z$-adequate. 
As a $\Z$-sup-preserving subset of $P$, $A$ already satisfies ${ x \ll^{\Z}_P y } \Rightarrow { x \ll^{\Z}_A y }$ for all $x, y \in A$ by Lemma~\ref{lem:reflecting}. 
To prove the reverse implication, let $x, y \in A$ such that $x \ll^{\Z}_A y$, and let $Z$ be a $\Z$-subset of $P$ with sup such that $y \leqslant \bigvee_P Z$. 
By the $\Z$-refinement property there is some $\Z$-subset $Z'$ of $A$ with sup in $A$ such that $Z' \subseteq { \downarrow_{P}\!\! Z }$ and $y \leqslant \bigvee_A Z'$. 
Since $x \ll^{\Z}_A y$, there is some $z' \in Z'$ such that $x \leqslant z'$. 
Moreover, $z' \in { \downarrow_{P}\!\! Z }$, so $x \leqslant z$ for some $z \in Z$. 
This entails $x \ll^{\Z}_P y$. 
So $A$ is a $\Z$-subposet of $P$. 

Now suppose that $A$ is a $\Z$-continuous $\Z$-subposet of $P$. 
Consider the inclusion map $i : A \to P$. 
This map is quasi-invertible with quasi-inverse the identity map of $A$. 
Since $A$ is a $\Z$-continuous $\Z$-subposet of $P$, the map $i$ has the $\Z$-refinement property by Proposition~\ref{prop:refinement}, which means that $A$ has the $\Z$-refinement property by Definition~\ref{def:refprop}.  
So $A$ is $\Z$-adequate. 
Suppose moreover that $\Z$ is union-complete and that $A$ has a weak $\Z$-basis $B$ and $\Z$-sup-generates $P$; then $P$ is $\Z$-continuous with weak $\Z$-basis $B$ by Theorem~\ref{thm:imageofzcp}. 
\end{proof}

\begin{theorem}\label{thm:union}
Let $(A_j)_{j \in J}$ be a covering family of subsets of a poset $P$, where each $A_j$ is a $\Z$-continuous $\Z$-subposet with weak $\Z$-basis $B_j$. 
Then $P$ is a $\Z$-continuous poset with weak $\Z$-basis $\bigcup_{j \in J} B_j$. 
\end{theorem}

\begin{proof}
Let $x \in P$. 
Let $k \in J$ such that $x \in A_{k}$. 
Since $A_k$ is a $\Z$-continuous poset with weak $\Z$-basis $B_k$, there is some $\Z$-subset $Z$ of $A_k$ included in 
\[
{ \downarrow_{A_k}\!\!(\twoheaddownarrow^{\Z}_{A_k} x \cap B_k) } = { \downarrow_{P}\!\! (\twoheaddownarrow^{\Z}_P x \cap B_k) \cap A_k }
\] 
whose sup in $A_k$ is $x$. 
Take $B := { \bigcup_{j \in J} B_j }$. 
Then $Z \subseteq { \downarrow_{P}\!\! (\twoheaddownarrow^{\Z}_P x \cap B) }$, and $Z$ is also a $\Z$-subset of $P$. 
%Since $C_k$ is included in $\downarrow\!\!(\twoheaddownarrow^{\Z}_P x \cap B)$, we have that $Z$ is also a $\Z$-subset of $\downarrow\!\!(\twoheaddownarrow^{\Z}_P x \cap B)$. 
Moreover, $A_k$ is $\Z$-sup-preserving in $P$, so $Z$ has a sup in $P$ equal to $x$. 
Thus, $P$ is $\Z$-continuous with weak $\Z$-basis $B$. 
\end{proof}

For the following corollary, recall that the disjoint union of a family $(P_j)_{j \in J}$ of posets is defined as the poset $\bigcup_{j \in J} \{ j \} \times P_j$ endowed with the partial order defined by $(j, x) \leqslant (k, y)$ if $j = k$ and $x \leqslant_{P_j} y$. 

\begin{corollary}
Let $(P_j)_{j \in J}$ be a family of posets, where each $P_j$ is $\Z^*$-continuous with weak $\Z^*$-basis $B_j$. 
Then the disjoint union of $(P_j)_{j \in J}$ is a $\Z^*$-continuous poset with weak $\Z^*$-basis the disjoint union of $(B_j)_{j \in J}$. 
\end{corollary}

\begin{proof}
Let $P$ be the disjoint union of $(P_j)_{j \in J}$, and take $A_j := \{ j \} \times P_j$. 
Then $P_j$ and $A_j$ are isomorphic posets, and $\{ j \} \times B_j$ is a weak $\Z^*$-basis of $A_j$, for all $j \in J$. 
To apply the previous theorem, we just need to show that every $A_j$ is a $\Z^*$-subposet of $P$. 
By Theorem~\ref{thm:main0}, it suffices to show that every $A_j$ is $\Z^*$-adequate in $P$. 
The latter holds by Example~\ref{ex:adequate}, since every $A_j$ is both a lower set and an upper set in $P$. 
\end{proof}

\subsection{Main theorem}

%\begin{definition}
%Given subsets $A$, $B$, and $C$ of a poset $P$ with $C \subseteq { A }$, we say that $(A,B)$ has the \textit{order-density property with respect to $C$} if $[a,b] \cap C$ is nonempty for all $a \in A$ and $b \in { B \cap \downarrow_{P}\!\! A }$ with $a \leqslant b$. 
%\end{definition}

%\begin{definition}
%Given subsets $A$, $B$, and $C$ of a poset $P$ with $C \subseteq { A }$, we say that $(A,B)$ has the \textit{order-density property with respect to $C$} if, whenever $a \leqslant b \ll^{\Z}_P a'$ for $a, a' \in A$ and $b \in B$, there exists some $c \in C$ with $a \leqslant c \ll^{\Z}_P a'$. 
%\end{definition}

The following definition can be seen as an extension of the notion of (strongly) $\Z$-dense subset (see Definition~\ref{def:dense}).

\begin{definition}
Given subsets $A$, $B$ of a poset $P$, a subset $C$ of $A$ is \textit{(strongly) $\Z$-dense in $A$ relatively to $(P, B)$}, or \textit{(strongly) $\Z$-dense in $A$ rel.\ $(P, B)$} for short, if, for all $a_1, a_2 \in A$ and $b \in B$,
\[
{ a_1 \prec b \ll^{\Z}_P a_2 } \Rightarrow { \exists c \in C, a_1 \prec c \ll^{\Z}_P a_2 }, 
\]
where $\prec$ is $\leqslant$ (resp.\ where $\prec$ is one of $\leqslant$ or $\ll^{\Z}_P$). 
%and \textit{strongly $\Z$-dense in $A \mid_{P} B$}, if it is $\Z$-dense in $A \mid_{P} B$ and, for all $a_1, a_2 \in A$ and $b \in B$,
%\[
%{ a_1 \ll^{\Z} b \ll^{\Z}_P a_2 } \Rightarrow { \exists c \in C, a_1 \ll^{\Z} c \ll^{\Z}_P a_2 }.
%\]
\end{definition}

%Be aware that a strongly $\Z$-dense subset in $A \mid_{P} B$ is not necessarily $\Z$-dense subset in $A \mid_{P} B$; yet we stick to the term \textit{strongly} for consistency with the notion of (strongly) $\Z$-dense subset of Definition~\ref{def:dense}. 

%\begin{example}\label{ex:odp}
%Let $A$ and $B$ be subsets of a poset $P$. 
%\begin{itemize}
%  \item $(A,B)$ has the order-density property with respect to $A \cap \downarrow_{P}\!\! B$. 
%  \item If $A$ is order-convex, then $(A, B)$ has the order-density property with respect to $A \cap B$. 
%  \item If $A$ is $\Z$-projected, and $p : P \to P$ is a projection that preserves the $\Z$-below relation such that $A = p(P)$, then $(A, B)$ has the order-density property with respect to $p(B)$. 
%  \item If $A$ is the direct image of a $\Z$-continuous poset by a $\Z$-adequate map, then $(A, B)$ has the order-density property with respect to $A \cap \uparrow_{P}\!\! B$. 
%  \item If $A$ is of the form $f(P')$, where $f : P' \to P$ is an order-preserving map such that there exists an order-preserving map $g : P \to P'$ with $f \circ g \circ f = f$ and $f(g(x)) \leqslant x$ for all $x \in P$, then $(A, B)$ has the order-density property with respect to $f(g(B))$. Indeed, let $a \in A$ and $b \in B$ with $a \leqslant b$. We can write $a = f(x')$ for some $x' \in P'$. So $a = f(x') = f(g(f(x'))) = f(g(a)) \leqslant f(g(b)) \leqslant b$, which implies the assertion.
%\end{itemize}
%\end{example}

A subset $U$ of a poset $P$ is \textit{$\Z$-Scott-open} if it is an upper set and $Z \cap U$ is nonempty whenever $Z$ is a $\Z$-subset of $P$ with sup such that $\bigvee_P Z \in U$. 
A subset of a poset $P$ is \textit{$\Z$-projected} if it is of the form $p(P)$, for some projection $p : P \to P$ preserving the $\Z$-below relation. 

\begin{theorem}[Main Theorem]\label{thm:main}
Let $A$ be a $\Z$-adequate subset of a poset $P$. 
Then $A$ is a $\Z$-subposet of $P$ and the $\Z$-Scott-open subsets of $A$ coincide with the subsets induced by the $\Z$-Scott-open subsets of $P$. 
Moreover, 
\begin{enumerate}
  \item\label{thm:main1} if $P$ is $\Z$-interpolating with $\Z$-dense subset $B$ and $C$ is strongly $\Z$-dense in $A$ rel.\ $(P, B)$, then $A$ is $\Z$-interpolating with $\Z$-dense subset $C$; 
  \item\label{thm:main3} if $P$ is (strongly) $\Z$-continuous with weak $\Z$-basis $B$ and $C$ is $\Z$-dense in $A$ rel.\ $(P, B)$, then $A$ is (strongly) $\Z$-continuous with weak $\Z$-basis $C\cap\uparrow_{P}\!\!B$; 
  \item\label{thm:main3bis} if $P$ is (strongly) $\Z$-continuous with weak $\Z$-basis $B$, then $A$ is (strongly) $\Z$-continuous with weak $\Z$-basis $A\cap\uparrow_{P}\!\!B$; 
  \item\label{thm:main2} if $P$ is (strongly) $\Z$-continuous, then $A$ is (strongly) $\Z$-continuous; 
  \item\label{thm:main4} if $P$ is $\Z$-algebraic and $A$ is either order-convex or $\Z$-projected, then $A$ is $\Z$-algebraic. 
\end{enumerate}
\end{theorem}

\begin{proof}
We already know from Theorem~\ref{thm:main0} that $A$ is a $\Z$-subposet of $P$. 
%As a $\Z$-sup-preserving subset of $P$, $A$ already satisfies $x \ll^{\Z}_P y \Rightarrow x \ll^{\Z}_A y$ for all $x, y \in A$ by Lemma~\ref{lem:reflecting}. 
%Conversely, let $x, y \in A$ such that $x \ll^{\Z}_A y$, and let $Z$ be a $\Z$-subset of $P$ with sup such that $y \leqslant \bigvee_P Z$. 
%By the $\Z$-refinement property there is some $\Z$-subset $Z'$ of $A$ with sup in $A$ such that $Z' \subseteq { \downarrow\!\! Z }$ and $y \leqslant \bigvee_A Z'$. 
%Since $x \ll^{\Z}_A y$, there is some $z' \in Z'$ such that $x \leqslant z'$. 
%Moreover, $z' \in { \downarrow\!\! Z }$, so $x \leqslant z$ for some $z \in Z$. 
%This entails $x \ll^{\Z}_P y$. 
%So $A$ is a $\Z$-subposet of $P$. 

Let $U$ be a $\Z$-Scott-open subset of $A$, and take $V := { \uparrow_{P}\!\! U }$. 
We prove that $V$ is $\Z$-Scott-open in $P$ and $U = { V \cap A }$. 
First, $V$ is an upper subset of $P$. 
Second, suppose that $\bigvee_P Z \in V$ for some $\Z$-subset $Z$ of $P$. 
Then $u \leqslant \bigvee_P Z$ for some $u \in U$. 
By the $\Z$-refinement property there is some $\Z$-subset $Z'$ of $A$ included in $\downarrow_{P}\!\! Z$, with sup in $A$ such that $u \leqslant \bigvee_A Z'$. 
Then $\bigvee_A Z' \in { \uparrow_A\!\! U } = { U }$; since $U$ is $\Z$-Scott-open in $A$ the intersection $Z' \cap U$ is nonempty. 
So let $z' \in { Z' \cap U }$. 
Then $z' \in { \downarrow_{P}\!\! Z }$, so there is some $z \in Z$ such that $z' \leqslant z$. 
This shows that $z \in { Z \cap \uparrow_{P}\!\! U } = { Z \cap V }$. 
So $Z \cap V$ is nonempty. 
This proves that $V$ is a $\Z$-Scott-open subset of $P$. 
Third, $U = { V \cap A }$, because ${ V \cap A } = { \uparrow_{P}\!\! U \cap A } = { \uparrow_A\!\! U } = { U }$. 

Let $V$ be a $\Z$-Scott-open subset of $P$, and let $U = { V \cap A }$. 
If $x \in { \uparrow_A\!\! U }$, then $u \leqslant x$ for some $u \in U$, so $u \in V$ and $x \in { \uparrow_{P}\!\! V } = { V }$, hence $x \in { V \cap A } = { U }$. 
So $U$ is an upper subset of $A$. 
If $\bigvee_A Z \in U$ for some $\Z$-subset $Z$ of $A$ with sup in $A$, then $\bigvee_P Z \in U$ since $A$ is $\Z$-sup-preserving, so $\bigvee_P Z \in V$. 
Since $V$ is $\Z$-Scott-open in $P$, the subset $Z \cap V$ is nonempty, so there is some $z \in { Z \cap V }$. 
So $z \in { V \cap A } = { U }$. 
This proves that $U$ is a $\Z$-Scott-open subset of $A$. 

Case \eqref{thm:main1}: 
Let $a_1, a_2 \in A$ such that $a_1 \ll^{\Z}_A a_2$. 
Since $A$ is a $\Z$-subposet, $a_1 \ll^{\Z}_P a_2$. 
Now $P$ is $\Z$-interpolating with $\Z$-dense subset $B$, so there is some $b \in B$ with $a_1 \ll^{\Z}_P b \ll^{\Z}_P a_2$. 
Using that $C$ is strongly $\Z$-dense in $A$ rel.\ $(P, B)$, there is some $c \in C$ with $a_1 \ll^{\Z}_P c \ll^{\Z}_P a_2$. 
Again, $A$ is a $\Z$-subposet, so $a_1 \ll^{\Z}_A c \ll^{\Z}_A a_2$. 
This proves that $A$ is $\Z$-interpolating with $\Z$-dense subset $C$. 

Case \eqref{thm:main3}: 
Let $x \in A$. 
Let $Z \subseteq { \downarrow_{P}\!\!(\twoheaddownarrow^{\Z}_P x \cap B) }$ be a $\Z$-subset of $P$ whose sup is $x$. 
By the $\Z$-refinement property there is some $\Z$-subset $Z'$ of $A$ included in $\downarrow_{P}\!\! Z$ with a sup in $A$ such that $x \leqslant \bigvee_A Z'$. 
Moreover, $Z' \subseteq { \downarrow_{P}\!\! Z } \subseteq { \downarrow_{P}\!\! x }$ and $x \in A$ imply $\bigvee_A Z' \leqslant x$, so that $x = \bigvee_A Z'$. 
Furthermore, $Z' \subseteq { \downarrow_{P}\!\!(\twoheaddownarrow^{\Z}_P x \cap B) \cap A }$. 
If $a \in { \downarrow_{P}\!\!(\twoheaddownarrow^{\Z}_P x \cap B) \cap A }$, then $a \leqslant b \ll^{\Z}_P x$ for some $b \in B$. 
By $\Z$-density of $C$ in $A$ rel.\ $(P, B)$, there is some $c \in C$ such that $a \leqslant c \ll^{\Z}_P x$. 
Thus, $c \ll^{\Z}_A x$. 
This shows that $a \in { \downarrow_A\!\!(\twoheaddownarrow^{\Z}_A x \cap C) }$. 
Consequently, $Z'$ is a $\Z$-subset of $A$ included in $\downarrow_A\!\!(\twoheaddownarrow^{\Z}_A x \cap C)$ whose sup in $A$ is $x$. 
This shows that $A$ is $\Z$-continuous with weak $\Z$-basis $C$. 

Suppose moreover that $P$ is $\Z$-interpolating, and let us show that $A$ is $\Z$-interpolating.  
Using Case \eqref{thm:main1}, it suffices to show that $C \cap \uparrow_{P}\!\! B$ is strongly $\Z$-dense in $A$ rel.\ $(P, B)$. 
So let $a_1, a_2 \in A$ and $b_1 \in B$ such that $a_1 \prec b_1 \ll^{\Z}_{P} a_2$, where ${ \prec } \in { \{ \leqslant, \ll^{\Z}_{P} \} }$. 
We already know that $A$ is $\Z$-continuous with weak $\Z$-basis $C$.  
Thus, we can write $a_2 = \bigvee_{A} Z_1$, for some $\Z$-subset $Z_1$ of $A$ included in $\downarrow_{A}\!\!(\twoheaddownarrow^{\Z}_{A} a_2 \cap C)$. 
Now $A$ is $\Z$-sup-preserving in $P$, so $a_2 = \bigvee_{P} Z_1$. 
Using $b_1 \ll^{\Z}_{P} a_2$, we get $b_1 \leqslant z$, for some $z \in Z_1$. 
From ${ Z_1 } \subseteq { \downarrow_{A}\!\!(\twoheaddownarrow^{\Z}_{A} a_2 \cap C) }$, we deduce that $a_1 \prec b_1 \leqslant z \leqslant c_1 \ll^{\Z}_{A} a_2$, for some $c_1 \in C$. 
Besides, $A$ is a $\Z$-subposet of $P$, so that $a_1 \prec c_1 \ll^{\Z}_{P} a_2$ and $c_1 \in { C \cap \uparrow_{P}\!\! B }$. 
This proves that $C \cap \uparrow_{P}\!\! B$ is strongly $\Z$-dense in $A$ rel.\ $(P, B)$, as required. 

Case \eqref{thm:main3bis} is a consequence of Case \eqref{thm:main3}, using that $A$ is $\Z$-dense in $A$ rel.\ $(P, B)$.  

Case \eqref{thm:main2} is a consequence of Case \eqref{thm:main3bis}: if $P$ is (strongly) $\Z$-continuous, then $P$ is (strongly) $\Z$-continuous with weak $\Z$-basis $P$ and $A$ is $\Z$-dense in $A$ rel.\ $(P, P)$, so $A$ is (strongly) $\Z$-continuous. 

Case \eqref{thm:main4}: 
By Case \eqref{thm:main2}, $A$ is $\Z$-continuous. 
Let us show that, whenever $x, y \in A$, $x \ll^{\Z}_A y$ implies $x \leqslant c \ll^{\Z}_{P} y$ for some $\Z$-compact element $c$ of $A$. 
Then we will have $c \ll^{\Z}_{A} y$ and, by Proposition~\ref{prop:basis}, this will be sufficient to prove that $A$ is $\Z$-algebraic.
So let $x, y \in A$ with $x \ll^{\Z}_A y$. 
Since $A$ is a $\Z$-subposet of $P$, we have $x \ll^{\Z}_P y$. 
By hypothesis, $P$ is $\Z$-algebraic, so there is some $\Z$-compact element $b$ of $P$ such that $x \leqslant b \ll^{\Z}_P y$. 
If $A$ is order-convex, then $b \in { [x, y]_P } \subseteq { A }$, so $b$ is $\Z$-compact in $A$. 
If $A$ is $\Z$-projected, we write $A = p(P)$ for some projection $p : P \to P$ preserving the $\Z$-below relation, so $x = p(x) \leqslant p(b) \ll^{\Z}_P p(y) = y$, and $p(b)$ is a $\Z$-compact element of $A$. 
\end{proof}

\begin{proposition}\label{prop:zdense}
Let $A$ and $B$ be subsets of a poset $P$. %, and let $C \subseteq A$.
%Suppose that one of the following conditions holds: 
\begin{enumerate}
  \item\label{prop:zdense1} $A \cap \downarrow_{P}\!\! B$ is $\Z$-dense in $A$ rel.\ $(P, B)$. 
  \item\label{prop:zdense6} If $A$ is of the form $f(P')$, where $f : P' \to P$ is the left adjoint of a pre-Galois connection $(f, g)$, then $f(g(B))$ is $\Z$-dense in $A$ rel.\ $(P, B)$. 
  \item\label{prop:zdense2} If $A$ is order-convex, then $A \cap B$ is strongly $\Z$-dense in $A$ rel.\ $(P, B)$. 
  \item\label{prop:zdense3} If $A$ is $\Z$-projected, and $p : P \to P$ is a projection that preserves the $\Z$-below relation such that $A = p(P)$, then $p(B)$ is strongly $\Z$-dense in $A$ rel.\ $(P, B)$. 
  \item\label{prop:zdense4} If $A$ $\Z$-sup-generates $P$ and $P$ is $\Z$-interpolating, then $A \cap \downarrow_{P}\!\!B$ is strongly $\Z$-dense in $A$ rel.\ $(P, B)$. 
%  \item If $A$ is a $\Z$-subposet of $P$ and $B$ is the set of $\Z$-compact elements of $P$, then $(A, B)$ has the order-density property with respect to the set of $\Z$-compact elements of $A$. Indeed,  let $a \in A$ and $b \in B \cap \downarrow\!\! A$ with $a \leqslant b$. Let $a' \in A$ such that $b \leqslant a'$. Since $b$ is $\Z$-compact in $P$, we have $a \leqslant b \ll^{\Z}_P b \leqslant a'$. Now $A$ is supposed to be a $\Z$-subposet of $P$, 
  \item\label{prop:zdense5} If $A$ is the direct image of a $\Z$-continuous poset by a $\Z$-adequate map, then $A \cap \uparrow_{P}\!\! B$ is strongly $\Z$-dense in $A$ rel.\ $(P, B)$. 
\end{enumerate}
%Then $C$ is $\Z$-dense in $A$ rel.\ $(P, B)$.
\end{proposition}

\begin{proof}
Case \eqref{prop:zdense1} is straightforward from the definitions. 

Case \eqref{prop:zdense6}:
Let $a_1, a_2 \in A$ and $b \in B$ with $a_1 \leqslant b \ll^{\Z}_{P} a_2$. 
We can write $a_1 = f(x')$ for some $x' \in P'$. 
Then $a_1 = f(x') = f(g(f(x'))) = f(g(a_1)) \leqslant f(g(b)) \leqslant b \ll^{\Z}_{P} a_2$, 
so that $a_1 \leqslant f(g(b)) \ll^{\Z}_{P} a_2$. 
This proves that $f(g(B))$ is $\Z$-dense in $A$ rel.\ $(P, B)$. 

For Cases \eqref{prop:zdense2}, \eqref{prop:zdense3}, \eqref{prop:zdense4}, and \eqref{prop:zdense5}, let ${ \prec } \in { \{ \leqslant, \ll^{\Z}_{P} \} }$, and let $a_1, a_2 \in A$ and $b \in B$ with $a_1 \prec b \ll^{\Z}_{P} a_2$. 

Case \eqref{prop:zdense2}:
From $a_1 \prec b \ll^{\Z}_{P} a_2$ we deduce that $a_1 \leqslant b \leqslant a_2$. 
Since $A$ is order-convex, we have $b \in { [a_1, a_2]_P } \subseteq { A }$. 
This proves that $A \cap B$ is strongly $\Z$-dense in $A$ rel.\ $(P, B)$.

Case \eqref{prop:zdense3}:
The subset $A$ is of the form $p(P)$ for some projection $p : P \to P$ preserving the $\Z$-below relation, so $a_1 = p(a_1) \prec p(b) \ll^{\Z}_P p(a_2) = a_2$. 
This proves that $p(B)$ is strongly $\Z$-dense in $A$ rel.\ $(P, B)$.

Case \eqref{prop:zdense4}:
If $\prec$ equals $\leqslant$, the result is clear by Case \eqref{prop:zdense1}. 
Now consider the case where $\prec$ equals $\ll^{\Z}_{P}$. 
Since $A$ $\Z$-sup-generates $P$, we can write $b = \bigvee_P Z$ for some $\Z$-subset $Z$ of $P$ included in $A$. 
By hypothesis, $P$ is $\Z$-interpolating, so there is some $x \in P$ such that $a_1 \ll^{\Z}_{P} x \ll^{\Z}_P b = \bigvee_P Z$. 
Thus, there is some $a \in Z$ such that $x \leqslant a$, so $a_1 \ll^{\Z}_{P} a \leqslant b \ll^{\Z}_P a_2$, which yields $a_1 \ll^{\Z}_{P} a \ll^{\Z}_P a_2$ and $a \in A \cap \downarrow_{P}\!\! B$. 
This proves that $A \cap \downarrow_{P}\!\! B$ is strongly $\Z$-dense in $A$ rel.\ $(P, B)$.

Case \eqref{prop:zdense5}:
The subset $A$ is of the form $f(P')$ for some $\Z$-continuous poset $P'$ and some $\Z$-adequate map $f : P' \to P$. 
We write $a_1 = f(x_1')$ and $a_2 = f(x_2')$, for some $x_1', x_2' \in P'$. 
Since $P'$ is $\Z$-continuous, we can write $x_2' = \bigvee_{P'} Z'$, for some $\Z$-subset $Z'$ of $P'$ included in $\twoheaddownarrow^{\Z}_{P'} x_2'$. 
Now $f$ is $\Z$-sup-preserving, so 
$
b \ll^{\Z}_{P} f(\bigvee_{P'} Z') = \bigvee_P f(Z')
$. 
So there is some $x' \in Z'$ with $b \leqslant f(x')$, and $x' \ll^{\Z}_{P'} x_2'$. 
Since $f$ preserves the $\Z$-below relation by Proposition~\ref{prop:refinement}, we have $f(x') \ll^{\Z}_P f(x_2') = a_2$. 
Moreover, $a_1 \prec b \leqslant f(x')$ implies that $a_1 \prec f(x')$. 
So $a_1 \prec f(x') \ll^{\Z}_P a_2$. 
Note also that $f(x') \in A \cap \uparrow_{P}\!\! B$. 
This proves that $A \cap \uparrow_{P}\!\! B$ is strongly $\Z$-dense in $A$ rel.\ $(P, B)$.
\end{proof}

%\begin{theorem}[Main Theorem, Part 1]\label{thm:mainpart1}
%Let $P$ be a $\Z$-interpolating poset, and $A$ be a $\Z$-subposet of $P$ satisfying one of the following conditions:
%\begin{enumerate}
%  \item\label{thm:mainpart11} $A$ is order-convex;
%  \item\label{thm:mainpart12} $A$ is $\Z$-sup-generating;
%  \item\label{thm:mainpart13} $A$ is $\Z$-projected;
%  \item\label{thm:mainpart14} $A$ is the direct image of a $\Z$-continuous poset by a $\Z$-adequate map. 
%\end{enumerate}
%Then $A$ is $\Z$-interpolating.
%\end{theorem}

In order to apply our Main Theorem, it becomes crucial to understand what types of subsets are $\Z$-adequate. 
We consider the case of lower subsets in Section~\ref{sec:lower}, order-convex subsets and $\Z$-Scott-convex subsets in Section~\ref{sec:convex}, upper subsets and $\D$-Scott open subsets in Section~\ref{sec:upper}, direct images of order-preserving maps in Section~\ref{sec:left}, kernel retracts in Section~\ref{sec:kernels}.

%%%%%%%%%%%%%%%%%%%%%%
%%%%%%%%%%%%%%%%%%%%%%
%%%%%%%%%%%%%%%%%%%%%%
%%%%%%%%%%%%%%%%%%%%%%
\section{Lower subsets as $\Z$-subposets}\label{sec:lower}

\subsection{$\Z$-lower-united subsets in arbitrary posets}

Let $P$ be a poset. 
Recall that a subset $L$ of $P$ is a lower set if $L = { \downarrow\!\! L }$, in other words, if $y \leqslant x \in L$ implies $y \in L$. 

\begin{lemma}\label{lem:supsemireg}
Every subset $L$ of a poset $P$ such that $P \!\setminus\! L \subseteq { L^{\uparrow} }$ is a lower set and is $\A$-sup-preserving. 
\end{lemma}

\begin{proof}
Let $x \in { \downarrow\!\! L }$, and suppose that $x \notin L$. 
Then there exists $\ell \in L$ such that $x \leqslant \ell$, and by assumption on $L$ we have $x \in L^{\uparrow}$. 
So $x = \ell \in L$, a contradiction. 
This proves that $x \in L$. 
So $L$ is a lower subset of $P$. 

Now let $A$ be a subset of $L$ with sup $\ell_0$ in $L$. 
Then $\ell_0$ is an upper bound of $A$ in $P$. 
Let $u \in P$ be another upper bound of $A$. 
If $u \in L$, then $\ell_0 \leqslant u$ by definition of $\ell_0$. 
If $u \notin L$, then $u \in L^{\uparrow}$, so again $\ell_0 \leqslant u$. 
This proves that $\ell_0$ is the sup of $A$ in $P$. 
So $L$ is $\A$-sup-preserving. 
\end{proof}

\begin{remark}
A converse statement to the previous lemma holds in chains: if $P$ is a chain, then every lower subset $L$ of $P$ satisfies $P \!\setminus\! L \subseteq { L^{\uparrow} }$. 
\end{remark}

We call a subset $L$ of a poset $P$ \textit{$\Z$-lower-united}, if for every $\Z$-subset $Z$ of $P$ with sup, $\bigvee_P Z \in P \!\setminus\! L$ implies $L \subseteq { \downarrow\!\! Z }$. 

\begin{theorem}\label{thm:lowerunited}
Every $\Z$-lower-united subset $L$ of a poset $P$ is a lower set, and is $\A$-sup-preserving and $\Z$-adequate. % (hence satisfies the conditions of Theorem~\ref{thm:main}). 
\end{theorem}

\begin{proof}
We first show that $L$ satisfies $P \!\setminus\! L \subseteq { L^{\uparrow} }$. 
Then by Lemma~\ref{lem:supsemireg}, this will imply that $L$ is a lower set and is $\A$-sup-preserving. 
To see this, let $x \in P \!\setminus\! L$. 
The subset $Z = \{ x \}$ is a $\Z$-subset of $P$, and $x = \bigvee_P Z \in { P \setminus L }$, hence $L \subseteq { \downarrow\!\! Z } = { \downarrow\!\! x }$. 
This shows that $x \in L^{\uparrow}$. 

Let $\ell \leqslant \bigvee_P Z$ for some $\ell \in L$ and some $\Z$-subset $Z$ of $P$ with sup. 
Suppose first that $\bigvee_P Z \in L$. 
Since $L$ is a lower set, this implies that $Z \subseteq { L }$. 
Take $Z' := Z$. 
Then we have $Z' \subseteq { \downarrow\!\! Z }$ and $\ell \leqslant \bigvee_P Z = \bigvee_L Z = \bigvee_L Z'$. 
Moreover, $Z'$ is a $\Z$-subset of $L$ by Proposition~\ref{prop:stable}\eqref{prop:stable6}. 
Now suppose that $\bigvee_P Z \in P \setminus L$. 
Since $L$ is $\Z$-lower-united, we get $L \subseteq { \downarrow\!\! Z }$. 
If $Z_{\ell}$ denotes the singleton $\{ \ell \}$, then $Z_{\ell}$ is a $\Z$-subset of $L$, ${ Z_{\ell} } \subseteq { \downarrow\!\! Z }$, and $\ell \leqslant \bigvee_L Z_{\ell} = \ell$. 
This shows that $L$ has the $\Z$-refinement property. 
So $L$ is $\Z$-adequate. 
\end{proof}

\subsection{Lower subsets in $\Z$-lower-Scott hereditary posets}

In \cite{Lawson04}, Lawson and Xu defined a \textit{poset with a lower hereditary Scott topology} as a poset such that, for every Scott-closed subset $A$, the relative Scott topology on $A$ agrees with the Scott topology of the poset $A$ \cite[Definition~3.1]{Lawson04}. 
See also Mao and Xu \cite{Mao09}. 
We transpose this notion to the framework of proper subset systems as follows. 
A subset $L$ of a poset $P$ is \textit{$\Z$-Scott-closed} if it is a lower set and $Z \subseteq { L }$ implies $\bigvee_P Z \in L$ for all $\Z$-subsets $Z$ of $P$ with sup. 
The $\Z$-Scott-closed subsets are exactly the complements of the $\Z$-Scott-open subsets. 
The family of $\Z$-Scott-closed subsets is stable under arbitrary intersections, but not under finite unions in general. 
A poset $P$ is \textit{$\Z$-lower-Scott hereditary} if, for every $\Z$-Scott-closed subset $L$, the $\Z$-Scott-closed subsets in the poset $L$ coincide with the subsets induced by the $\Z$-Scott-closed subsets of $P$. 
The following extends Lawson and Xu \cite[Lemma~3.2]{Lawson04}. 

\begin{lemma}\label{lem:characzlsh}
Let $P$ be a poset. 
The following are equivalent: 
\begin{enumerate}
  \item\label{lem:characzlsh1} $P$ is $\Z$-lower-Scott hereditary; 
  \item\label{lem:characzlsh2} every principal ideal of $P$ is $\Z$-sup-preserving; 
  \item\label{lem:characzlsh3} every lower subset of $P$ is $\Z$-sup-preserving; 
  \item\label{lem:characzlsh4} any minimal upper bound of $Z$ is actually the sup of $Z$, for every $\Z$-subset $Z$ of $P$. 
\end{enumerate}
\end{lemma}

\begin{proof}
\eqref{lem:characzlsh1} $\Rightarrow$ \eqref{lem:characzlsh2}: 
Let $x \in P$ and take $L := { \downarrow_{P}\!\! x }$. 
Then $L$ is a $\Z$-Scott-closed subset of $P$. 
Let $Z$ be a $\Z$-subset of $L$ with sup $z_0$ in $L$, and let $u$ be an upper bound of $Z$ in $P$. 
Then $\downarrow_{P}\!\! u \cap L$ is $\Z$-Scott-closed in the poset $L$ and $Z \subseteq { \downarrow_{P}\!\! u \cap L }$, so $z_0 \in { \downarrow_{P}\!\! u \cap L }$. 
In particular, $z_0 \leqslant u$, so $z_0$ is the sup of $Z$ in $P$. 
Thus, $L$ is $\Z$-sup-preserving. 

\eqref{lem:characzlsh2} $\Rightarrow$ \eqref{lem:characzlsh4}: 
Let $u$ be a minimal upper bound of a $\Z$-subset $Z$ of $P$. 
Then $u$ is the sup of $Z$ in $\downarrow_{P}\!\! u$. Since $\downarrow_{P}\!\! u$ is $\Z$-sup-preserving by hypothesis, $u$ is also the sup of $Z$ in $P$. 

\eqref{lem:characzlsh4} $\Rightarrow$ \eqref{lem:characzlsh3}: 
Let $L$ be a lower subset of $P$, and let $Z$ be a $\Z$-subset of $L$ with sup $\ell$ in $L$. 
We show that $\ell$ is a minimal upper bound of $Z$ in $P$. 
So let $u$ be an upper bound of $Z$ in $P$ such that $u \leqslant \ell$. 
Since $L$ is a lower set, we have $u \in L$. 
By definition of $\ell$, this implies that $\ell \leqslant u$, hence $u = \ell$. 
So $\ell$ is a minimal upper bound of $Z$ in $P$, hence is the sup of $Z$ in $P$. 

\eqref{lem:characzlsh3} $\Rightarrow$ \eqref{lem:characzlsh1}: 
Let $L$ be a $\Z$-Scott-closed subset of $P$. 
First we consider a $\Z$-Scott-closed subset $A$ of $L$, and we show that $A$ is a $\Z$-Scott-closed subset of $P$. 
We have $A \subseteq { \downarrow_{P}\!\! A } \subseteq { \downarrow_{P}\!\! L \cap \downarrow_{P}\!\! A } = { L \cap \downarrow_{P}\!\! A } = { \downarrow_{L}\!\! A } = { A }$, so $A$ is a lower subset of $P$. 
If $Z$ is a $\Z$-subset of $P$ with sup $z_0$ such that $Z \subseteq { A }$, then $z_0$ belongs to $L$ since $L$ is a $\Z$-Scott-closed subset of $P$. 
This implies that $z_0$ is also the sup of $Z$ in $L$; since $A$ is a $\Z$-Scott-closed subset of $L$ we obtain $z_0 \in A$. 
So $A$ is indeed a $\Z$-Scott-closed subset of $P$. 

Conversely, let $A$ be a subset of $P$ of the form $F \cap L$, where $F$ is a $\Z$-Scott-closed subset of $P$. 
We show that $A$ is a $\Z$-Scott-closed subset of $L$. 
We have $A \subseteq { \downarrow_{L}\!\! A } = { \downarrow_{P}\!\! A \cap L } \subseteq { \downarrow_{P}\!\! F \cap \downarrow_{P}\!\! L \cap L } = { F \cap L } = { A }$, so $A$ is a lower subset of $L$. 
If $Z$ is a $\Z$-subset of $L$ with sup $z_0 \in L$ such that $Z \subseteq { A }$, then $z_0 = \bigvee_P Z$ by \eqref{lem:characzlsh3}. 
This implies that $z_0 \in F$ since $F$ is a $\Z$-Scott-closed subset of $P$. 
So $z_0 \in F \cap L = A$. 
This shows that $A$ is a $\Z$-Scott-closed subset of $L$. 
\end{proof}

In addition to the previous lemma, we shall see (with Proposition~\ref{prop:zriesz}) sufficient conditions for a poset to be $\Z$-lower-Scott hereditary; in particular, we shall prove that semilattices and conditionally $\Z$-complete posets are $\Z$-lower-Scott hereditary. 

\begin{definition}
A subset $A$ of a poset $P$ is \textit{$\Z$-sup-regular} if, for every $\Z$-subset $Z$ of $A$, $Z$ has a sup in $A$ if and only if $Z$ has a sup in $P$,  and $\bigvee_P Z = \bigvee_A Z$ if one of these sups exists. 
\end{definition}

Every $\Z$-sup-regular subset is of course $\Z$-sup-preserving. 

%\begin{example}
%A subset $A$ of a poset $P$ is called an \textit{interior system} of $P$ if, for every $x \in P$, there exists a smallest $a \in A$ with $x \leqslant a$. 
%\end{example}

\begin{lemma}\label{lem:lower}
Every $\Z$-Scott-closed subset of a $\Z$-lower-Scott hereditary poset $P$ is $\Z$-sup-regular.  
\end{lemma}

\begin{proof}
Let $L$ be a $\Z$-Scott-closed subset of $P$. 
We know with Lemma~\ref{lem:characzlsh} that $L$, as a lower set, is $\Z$-sup-preserving. 
Now if $Z$ is a $\Z$-subset of $L$ with sup in $P$, then $\bigvee_P Z \in L$ by definition of $\Z$-Scott-closedness, so $\bigvee_P Z = \bigvee_L Z$. 
So $L$ is $\Z$-sup-regular. 
\end{proof}

The following result extends Ern\'e \cite[Proposition~1]{Erne05}. 

\begin{theorem}\label{thm:lower}
Every lower subset of a $\Z$-continuous, $\Z$-lower-Scott hereditary poset $P$ is $\Z$-adequate. 
\end{theorem}

\begin{proof}
Let $L$ be a lower subset of $P$. 
By Lemma~\ref{lem:characzlsh}, $L$ is $\Z$-sup-preserving. 
Let us show that $L$ has the $\Z$-refinement property. 
Let $\ell \in L$ and $Z$ be a $\Z$-subset of $P$ such that $\ell \leqslant \bigvee_P Z$. 
Then $\twoheaddownarrow^{\Z}_P \ell \subseteq { \downarrow\!\! Z \cap \downarrow\!\! \ell } \subseteq { \downarrow\!\! Z \cap L }$, using the definition of the $\Z$-below relation and the fact that $L$ is a lower set. 
Moreover, $P$ is $\Z$-continuous, so there is some $\Z$-subset $Z'$ of $P$ included in $\twoheaddownarrow^{\Z}_P \ell$ with sup $\ell$. 
So $Z' \subseteq { \downarrow\!\! Z \cap L }$ and $\ell = \bigvee_P Z' = \bigvee_L Z'$, 
Also, $Z'$ is a $\Z$-subset of $L$ by Proposition~\ref{prop:stable}\eqref{prop:stable6}. 
So $L$ has the $\Z$-refinement property. 
\end{proof}

We also state the following variation of the previous theorem.

\begin{theorem}
Every subset $L$ of a $\Z$-continuous poset $P$ with $P \setminus L \subseteq { L^{\uparrow} }$ is $\A$-sup-preserving and $\Z$-adequate. % (hence satisfies the conditions of Theorem~\ref{thm:main}). 
\end{theorem}

\begin{proof}
The subset $L$ is a lower set and is $\A$-sup-preserving by Lemma~\ref{lem:supsemireg}. 
One can then show that $L$ has the $\Z$-refinement property along the same lines as in the previous proof. 
\end{proof}

The following extends Lawson and Xu \cite[Lemma~3.4]{Lawson04}. 

\begin{corollary}\label{coro:lawsonxu}
Let $P$ be a $\Z$-lower-Scott hereditary poset. 
The following are equivalent: 
\begin{enumerate}
  \item\label{coro:lawsonxu1} every principal ideal of $P$ is a (strongly) $\Z$-continuous poset; 
  \item\label{coro:lawsonxu2} every lower subset of $P$ is a (strongly) $\Z$-continuous poset; 
  \item\label{coro:lawsonxu3} $P$ is a (strongly) $\Z$-continuous poset. 
\end{enumerate}
\end{corollary}

\begin{proof}
\eqref{coro:lawsonxu3} $\Rightarrow$ \eqref{coro:lawsonxu2} is a consequence of Theorem~\ref{thm:lower}, once combined with Theorem~\ref{thm:main}. 

\eqref{coro:lawsonxu2} $\Rightarrow$ \eqref{coro:lawsonxu1} is straightforward. 

\eqref{coro:lawsonxu1} $\Rightarrow$ \eqref{coro:lawsonxu3}:
Let $x \in P$. 
Then by \eqref{coro:lawsonxu1} the principal ideal $\downarrow\!\! x$ is a $\Z$-continuous poset. 
By Lemma~\ref{lem:characzlsh}, $\downarrow\!\! x$ is $\Z$-sup-preserving. 
Let us show that $\downarrow\!\! x$ has the $\Z$-refinement property. 
So let $y \in { \downarrow\!\! x }$, and let $Z$ be a $\Z$-subset of $P$ with sup $z$ such that $y \leqslant z = \bigvee_P Z$. 
Take $L := { \downarrow\!\! z }$. 
Then $Z \subseteq { L }$ and we have $z = \bigvee_L Z$. 
Moreover, $Z$ is a $\Z$-subset of $L$ by Proposition~\ref{prop:stable}\eqref{prop:stable6}. 
Now, $y \leqslant \bigvee_L Z$ implies $\twoheaddownarrow^{\Z}_L y \subseteq { \downarrow\!\! Z \cap L }$. 
Since $L$ is $\Z$-continuous, there is some $\Z$-subset $Z_y$ of $L$ such that $Z_y \subseteq { \twoheaddownarrow^{\Z}_L y }$ and $\bigvee_L Z_y = y$. 
We have $y \leqslant x$, so $Z_y \subseteq { \downarrow\!\! x }$. 
Also, $y \in Z_y^{\uparrow} \cap \downarrow\!\! x$, so $Z_y$ is a $\Z$-subset of $\downarrow\!\! x$ by Proposition~\ref{prop:stable}\eqref{prop:stable6}. 
By Lemma~\ref{lem:characzlsh}, $L$ is $\Z$-sup-preserving, so that $\bigvee_L Z_y = \bigvee_P Z_y = y$. 
Thus, $\bigvee_P Z_y = y \in { \downarrow\!\! x }$ and $Z_y \subseteq { \downarrow\!\! x }$ imply $y = \bigvee_{\downarrow x} Z_y$. 
This proves that $\downarrow\!\! x$ has the $\Z$-refinement property. 

By Theorem~\ref{thm:main}, we deduce that every principal ideal of $P$ is a $\Z$-continuous $\Z$-subposet of $P$. 
Now Theorem~\ref{thm:union} applies, and we conclude that $P$ is $\Z$-continuous. 

To finish the proof, suppose that every principal ideal of $P$ is a strongly $\Z$-continuous poset, 
and let $x, z \in P$ with $x \ll^{\Z}_P z$. 
Then $x, z \in L$ with ${ L } := { \downarrow\!\! z }$. 
From the above we know that $L$ is a $\Z$-subposet of $P$, so $x \ll^{\Z}_L z$. 
Since $L$ is $\Z$-interpolating by hypothesis, $x \ll^{\Z}_L y \ll^{\Z}_L z$ for some $y \in L$. 
So $x \ll^{\Z}_P y \ll^{\Z}_P z$, which proves that $P$ is $\Z$-interpolating, as required. 
\end{proof}

\subsection{$\Z$-Riesz posets}

A poset $P$ is a \textit{semilattice} if every nonempty finite subset has an infimum; we write $x \wedge y$ for the infimum of $\{ x, y \}$. 
A nonempty subset $A$ of a poset is \textit{filtered} if every finite subset of $A$ has a lower bound in $A$. 
A poset is \textit{$\Z$-Riesz} if $Z^{\uparrow}$ is filtered for every upper-bounded $\Z$-subset $Z$. 

\begin{proposition}\label{prop:zriesz}
The following assertions hold: 
\begin{enumerate}
  \item\label{prop:zriesz1} Every conditionally $\Z$-complete poset is $\Z$-Riesz. 
  \item\label{prop:zriesz2} Every semilattice is $\Z$-Riesz. 
  \item\label{prop:zriesz3} Every $\Z$-Riesz poset is $\Z$-lower-Scott hereditary. 
  \item\label{prop:zriesz4} A poset is $\A$-Riesz if and only if it is a semilattice. 
\end{enumerate}
\end{proposition}

\begin{proof}
Case \eqref{prop:zriesz1} is clear. 

Case \eqref{prop:zriesz2} is implied by Case \eqref{prop:zriesz4}. 

Case \eqref{prop:zriesz3} is a direct consequence of Lemma~\ref{lem:characzlsh}\eqref{lem:characzlsh4}. 

Case \eqref{prop:zriesz4}:
The fact that a semilattice is $\A$-Riesz (hence $\Z$-Riesz for every proper subset system $\Z$) is straightforward. 
Let us prove that an $\A$-Riesz poset is necessarily a semilattice. 
Let $F$ be a nonempty finite subset of an $\A$-Riesz poset $P$. 
We want to show that $F$ has an infimum. 
Let $A = F^{\downarrow}$. 
Then $F \subseteq { A^{\uparrow} }$, so $A^{\uparrow}$ is (nonempty) filtered. 
Thus, there exists some $z \in A^{\uparrow}$ such that $z \leqslant f$ for all $f \in F$. 
In particular, $z \in F^{\downarrow} = A$. 
Now $z \in A \cap A^{\uparrow}$ means that $z$ is the greatest element of $A$, i.e.\ the greatest lower bound of $F$, i.e.\ the infimum of $F$. 
\end{proof}

\begin{corollary}
Every lower subset of a semilattice is $\A$-sup-preserving. 
\end{corollary}

%%%%%%%%%%%%%%%%%%%%%%
%%%%%%%%%%%%%%%%%%%%%%
%%%%%%%%%%%%%%%%%%%%%%
%%%%%%%%%%%%%%%%%%%%%%
\section{Order-convex subsets as $\Z$-subposets in semilattices}\label{sec:convex}

Recall that a subset $C$ of a poset is order-convex if and only if $C$ is the intersection of a lower subset with an upper subset. 
%A poset $P$ is a \textit{semilattice} if every nonempty finite subset has an infimum; we write $x \wedge y$ for the infimum of $\{ x, y \}$. 

\begin{lemma}\label{lem:convex}
Every order-convex subset of a $\Z$-Riesz poset $P$ is $\Z^*$-sup-preserving. 
\end{lemma}

\begin{proof}
%Let $P$ be a $\Z$-Riesz poset and 
Let $C$ be an order-convex subset of $P$. 
Let $Z$ be a nonempty $\Z$-subset of $C$ with sup $c_0$ in $C$. 
Let $u \in P$ be an upper bound of $Z$ in $P$. 
Since $Z$ is nonempty there is some $z \in Z$. 
By hypothesis, $P$ is $\Z$-Riesz, so the subset $Z^{\uparrow}$ is filtered. 
Thus, there exists some $v \in Z^{\uparrow}$ such that $v \leqslant u$ and $v \leqslant c_0$. 
In particular, $z \leqslant v \leqslant c_0$, so $v \in C$ by order-convexity of $C$. 
This implies that $v$ is an upper bound of $Z$ in $C$, hence $c_0 \leqslant v \leqslant u$. 
This shows that $c_0$ is the sup of $Z$ in $P$. 
So $C$ is $\Z^*$-sup-preserving. 
\end{proof}

We say that a subset of a poset is \textit{$\Z$-Scott-convex} if it is the intersection of a lower subset with a $\Z$-Scott-open subset of the poset. 
Hence every $\Z$-Scott-convex subset of a poset is order-convex. 
A semilattice is \textit{$\Z$-meet-continuous} if 
\[
(\bigvee Z) \wedge x = \bigvee (Z \wedge x), 
\]
for every $x$ and $\Z$-subset $Z$ with sup, where $Z \wedge x$ denotes the subset 
\[
Z \wedge x := \{ z \wedge x : z \in Z \}.
\] 
Note that $Z \wedge x$ is necessarily a $\Z$-subset as the image of $Z$ by the order-preserving map $y \mapsto y \wedge x$. 
It is not difficult to show that every semilattice that is a $\Z$-continuous poset is $\Z$-meet-continuous (see e.g.\ \cite[Proposition~I-1.8]{Gierz03} for a proof in the frame of classical domain theory). 

\begin{theorem}\label{thm:scottconvex}
Every $\Z$-Scott-convex subset $C$ of a $\Z$-meet-continuous semilattice $P$ is order-convex, $\A^*$-sup-preserving, and $\Z$-adequate. %(hence satisfies the conditions of Theorem~\ref{thm:main} if $\Z = \A^*$). 
\end{theorem}

\begin{proof}
We can write $C = { \downarrow\!\! C \cap U }$ for some $\Z$-Scott-open subset $U$ of $P$, by definition of $\Z$-Scott-convexity. 
As a semilattice, $P$ is $\A$-Riesz by Proposition~\ref{prop:zriesz}, so the subset $C$ is $\A^*$-sup-preserving by Lemma~\ref{lem:convex}. 
Now, we show that $C$ has the $\Z$-refinement property. 
To see this, let $c \leqslant \bigvee_P Z$ for some $c \in C$ and some $\Z$-subset $Z$ of $P$ with sup. 
If $Z' := { Z \wedge c }$, then $Z'$ is a $\Z$-subset of $P$ as noticed in the lines preceding the theorem, $\bigvee_P Z' = c$ by $\Z$-meet-continuity of $P$, and $Z' \subseteq { \downarrow\!\! Z }$. 
In particular, $\bigvee_P Z' \in U$, so $Z' \cap U$ is nonempty since $U$ is $\Z$-Scott-open. 
We pick some $u \in { Z' \cap U }$; note that $u \leqslant c$. 
We now consider $Z'' := { (Z' \cap \uparrow\!\! u) \wedge c }$. 
Then $Z''$ is a $\Z$-subset of $C$. 
Indeed, let $f_u : P \to { \uparrow\!\! u }$ be the order-preserving map defined in the proof of Proposition~\ref{prop:stable}, and let $g : { \uparrow\!\! u } \to C$ be defined by $g(y) = y \wedge c$. 
Note that $g$ is well-defined, in the sense that $g(y) \in C$ for all $y \in { \uparrow\!\! u }$. 
This is because, if $y \in { \uparrow\!\! u }$, then $g(y) = y \wedge c \in { [u, c]_{P} }$; noting that $c \in C$, $u \in { \downarrow\!\! C \cap U } = { C }$, and $C$ is order-convex, we obtain $g(y) \in C$. 
Moreover, $g$ is obviously order-preserving. 
Then $Z'' = { g \circ f_u(Z') }$, so $Z''$ is a $\Z$-subset of $C$. 
%If $x \in Z''$, then $x \leqslant c$ on the one hand, and $u = u \wedge c \leqslant x$ on the other hand, so $x \in { \downarrow\!\! C \cap \uparrow\!\! U } = { \downarrow\!\! C \cap U } = { C }$. 
%This implies that $Z'' \subseteq { \downarrow\!\! Z \cap C }$. 
Besides, $Z'' \subseteq { \downarrow\!\! Z }$ and $c = \bigvee_C Z''$, so we have shown that $C$ has the $\Z$-refinement property. 
%
%Proving that, if $C$ is a subsemilattice of $P$, then $C$ is a $\Z$-meet-continuous semilattice, is straightforward and left to the reader. 
\end{proof}

\begin{remark}
If $C$ is a subsemilattice and a $\Z$-sup-preserving subset of a $\Z$-meet-continuous semilattice $P$, it is easy to show that $C$ is a $\Z$-meet-continuous semilattice itself. 
\end{remark}

\begin{corollary}
Every lower subset of a $\Z$-meet-continuous semilattice is order-convex, $\A$-sup-preserving, and $\Z$-adequate. % (hence satisfies the conditions of Theorem~\ref{thm:main}). 
\end{corollary}

%%%%%%%%%%%%%%%%%%%%%%
%%%%%%%%%%%%%%%%%%%%%%
%%%%%%%%%%%%%%%%%%%%%%
%%%%%%%%%%%%%%%%%%%%%%
\section{$\D$-Scott-open subsets as $\D$-subposets}\label{sec:upper}

In this section, we specialize to the case $\Z = \D$ and examine properties of upper subsets and $\D$-Scott-open subsets. 

\begin{lemma}\label{lem:sup}
Let $P$ be a poset and $U$ be an upper subset of $P$. 
If $D$ is a directed subset of $P$ with sup such that $D \cap U$ is nonempty, then $D \cap U$ is a directed subset of $P$ with sup and $\bigvee (D \cap U) = \bigvee D$. 
\end{lemma}

\begin{proof}
First, we show that the nonempty subset $D \cap U$ is directed. 
To see this, let $d, d' \in D \cap U$. 
Since $D$ is directed, there is some $d'' \in D$ such that $d \leqslant d''$ and $d' \leqslant d''$. 
Moreover, $d'' \in U$ for $U$ is an upper set. 
So $d'' \in D \cap U$. 
This shows  that $D \cap U$ is directed. 

Second, we prove that $\bigvee (D \cap U) = \bigvee D$. 
It is trivial that $\bigvee D$ is an upper bound of $D \cap U$. 
Further, let $c \in P$ be another upper bound of $D \cap U$. 
Since $D \cap U$ is nonempty by hypothesis, we pick some $d_0 \in D \cap U$. 
Now let $d \in D$. 
Then there is some $d_1 \in D$ such that $d \leqslant d_1$ and $d_0 \leqslant d_1$. 
Since $U$ is an upper set, $d_1 \in U$, so $d_1 \in D \cap U$. 
Thus, we have $d \leqslant d_1 \leqslant c$. 
So $d \leqslant c$, for all $d \in D$. 
This gives $\bigvee D \leqslant c$, as required. 
\end{proof}

%\begin{proof}
%Since $D \cap U$ is nonempty by hypothesis, we pick some $d_0 \in D \cap U$. 
%If $d, d' \in D$, then there exists some $d_1 \in D$ greater than $d_0$, $d$, and $d'$. 
%In particular, $d_1 \in D \cap \uparrow\!\!\{d_0\} \subseteq { D \cap \uparrow\!\! U } = { D \cap U }$ since $U$ is an upper set. 
%This shows that $D \cap U$ is directed and $D \subseteq { \downarrow\!\! (D \cap U) } \subseteq { \downarrow\!\! D }$, so that $\bigvee D = \bigvee (D \cap U)$. 
%\end{proof}

\begin{remark}
The property of the previous lemma actually characterizes upper subsets of a poset. 
Indeed, assume that this property is satisfied, and let $y \in { \uparrow\!\! U }$. 
Then $x \leqslant y$ for some $x \in U$. 
The subset $D = \{ x, y \}$ is directed, its sup is $y$, and $D \cap U$ contains $x$ so is nonempty. 
Thus, $\bigvee D = y = \bigvee (D \cap U)$. 
This latter equality implies that $y \in U$. 
\end{remark}

\begin{theorem}\label{thm:sup2}
Every upper subset of a poset is $\A^*$-sup-regular. 
Moreover, it is $\D$-Scott-open if and only if it has the $\D$-refinement property (in this case it satisfies the conditions of Theorem~\ref{thm:main} for $\Z = \D$). 
\end{theorem}

\begin{proof}
Let $U$ be an upper subset of a poset $P$, and let $A$ be a nonempty subset of $U$. 
Suppose first that $A$ has a sup $x_0$ in $U$. 
Then $x_0$ is an upper bound of $A$ in $P$. 
Let $u \in P$ be another upper bound of $A$ in $P$. 
Since $A$ is nonempty, this entails $u \in { \uparrow_{P}\!\! A } \subseteq { \uparrow_{P}\!\! U } = { U }$. 
Thus, $x_0 \leqslant u$, which shows that $x_0$ is the sup of $A$ in $P$, i.e.\ $\bigvee_U A = x_0 = \bigvee_P A$. 

Now suppose that $A$ has a sup $x_0$ in $P$. 
Again $A$ is nonempty, so $x_0 \in { \uparrow_{P}\!\! A } \subseteq { \uparrow_{P}\!\! U } = { U }$. 
This implies that $x_0$ is an upper bound of $A$ in $U$. 
Let $u \in U$ be another upper bound of $A$ in $U$. 
Then $u$ is an upper bound of $A$ (in $P$), so $x_0 \leqslant u$. 
This shows that $x_0$ is the sup of $A$ in $U$, i.e.\ $\bigvee_P A = x_0 = \bigvee_U A$.
At this stage we have proved that $U$ is $\A^*$-sup-regular. 

%As an upper subset, $U$ is regular by Lemma~\ref{lem:sup0}. 
Suppose that $U$ is $\D$-Scott-open. 
Let $u \leqslant \bigvee_P D$ for some $u \in U$ and some directed subset $D$ of $P$ with sup. 
Since $U$ is an upper set, $\bigvee_P D \in U$, and so $\bigvee_P D = \bigvee_U D$. 
Since $U$ is $\D$-Scott-open, the subset $D' := { D \cap U }$ is nonempty. 
By Lemma~\ref{lem:sup}, $D'$ is a directed subset of $P$ (hence of $U$) and $\bigvee_P D = \bigvee_P (D \cap U) = \bigvee_P D'$. 
Since $\bigvee_P D' = \bigvee_P D \in U$, we obtain $\bigvee_U D' = \bigvee_P D$. 
Thus, $u \leqslant \bigvee_P D = \bigvee_U D'$ and also we have $D' \subseteq { D } \subseteq { \downarrow_{P}\!\! D }$. 
This shows that $U$ has the $\D$-refinement property. 

Conversely, suppose that $U$ has the $\D$-refinement property. 
Let $D$ be a directed subset of $P$ with sup such that $\bigvee_P D \in U$. 
By the $\D$-refinement property there is a directed subset $D'$ of $U$ such that $D' \subseteq { \downarrow_{P}\!\! D }$. 
As a directed set, $D'$ is nonempty, so that ${ D \cap U } = { D \cap \uparrow_{P}\!\! U } \neq { \emptyset }$. 
Since $U$ is already an upper set by hypothesis, this proves that $U$ is $\D$-Scott-open. 
\end{proof}

\begin{example}\label{ex:scottopen}
Let $P$ be a $\D$-interpolating poset. 
Then every subset of the form 
\[
\twoheaduparrow{}^{\D} \! A := \{ x \in P : a \ll^{\D} x \mbox{ for some } a \in A \}, 
\]
with $A \subseteq { P }$, is $\D$-Scott-open, hence is $\A^*$-sup-regular and $\D$-adequate. 
\end{example}

\begin{corollary}
Let $P$ be a poset and $x \in P$. 
Then the principal filter $\uparrow\!\! x$ is a $\D$-subposet of $P$ if and only if $x$ is $\D$-compact. 
\end{corollary}

\begin{proof}
We already know that $\uparrow\!\! x$ is $\A^*$-sup-regular by Theorem~\ref{thm:sup2}. 
Moreover, using that $x$ is $\D$-compact it is easy to see that $\uparrow\!\! x$ is $\D$-Scott-open; so $\uparrow\!\! x$ has the $\D$-refinement property by Theorem~\ref{thm:sup2}, hence is $\D$-adequate. 
Then $\uparrow\!\! x$ is a $\D$-subposet of $P$ by Theorem~\ref{thm:main}. 

Conversely, suppose that $A := { \uparrow\!\! x }$ is a $\D$-subposet of $P$. 
As the least element of $A$, $x$ satisfies $x \ll^{\D}_A x$. 
Since $A$ is a $\D$-subposet of $P$, this implies $x \ll^{\D}_P x$, i.e.\ $x$ is $\D$-compact. 
\end{proof}

Putting Theorem~\ref{thm:scottconvex} and Theorem~\ref{thm:sup2} together, we obtain the following result. 

\begin{corollary}
Every $\Z$-Scott-open subset of a $\Z$-meet-continuous semilattice $P$ is $\A^*$-sup-regular and $\Z$-adequate. 
\end{corollary}

%%%%%%%%%%%%%%%%%%%%%%
%%%%%%%%%%%%%%%%%%%%%%
%%%%%%%%%%%%%%%%%%%%%%
%%%%%%%%%%%%%%%%%%%%%%
\section{Direct images of posets}\label{sec:left}

In this section, we study subsets of the form $f(P)$, for some order-preserving map $f : P \to P'$. 
Our interest is twofold: 
\begin{itemize}
  \item finding conditions under which $f(P)$ is $\Z$-adequate, in which case $f(P)$ becomes $\Z$-continuous as soon as $P'$ is $\Z$-continuous, by Theorem~\ref{thm:main} (Subsection~\ref{sec:left1});
  \item finding conditions under which $f(P)$ is $\Z$-continuous as soon as $P$ is $\Z$-continuous (Subsection~\ref{sec:left2}). 
\end{itemize}
We also examine properties of subsets of the form $f(A)$ with $A \subseteq P$ (Subsection~\ref{sec:left3}). 

\subsection{Direct images as $\Z$-subposets}\label{sec:left1}

The following result provides us with sufficient conditions for a subset of the form $f(P)$, where $f : P \to P'$ is an order-preserving map, to be $\Z$-adequate in $P'$. 

\begin{theorem}\label{thm:adequate}
Let $f : P \to P'$ be an order-preserving map. 
Assume that one of the following conditions holds: 
\begin{enumerate}
  \item\label{thm:adequate1} $f$ is surjective; 
  \item\label{thm:adequate2} $f$ is $\Z$-adequate and $f(P)$ is $\Z$-sup-preserving in $P'$;
  \item\label{thm:adequate3} $f$ is the left adjoint of a Galois connection $(f, g)$ where $g$ is $\Z$-sup-preserving;
  \item\label{thm:adequate4} $f$ is $\Z$-sup-preserving and is the left adjoint of a pre-Galois connection $(f, g)$ where $g$ is $\Z$-sup-preserving; 
  \item\label{thm:adequate6} $P$ is cond.\ $\Z$-complete and $f$ is a $\Z$-adequate projection;
  \item\label{thm:adequate5} $P$ is cond.\ $\Z$-complete and $f$ is $\Z$-adequate quasi-invertible;
  \item\label{thm:adequate7} $f$ is a $\Z$-adequate order-embedding.
%  \item\label{thm:adequate8} $P'$ is conditionally $\Z$-complete, $f$ is $\Z$-sup-preserving, and there is a $\Z$-sup-preserving map $g : P' \to P$ such that $f \circ g \circ f = f$. 
\end{enumerate} 
Then $f(P)$ is $\Z$-adequate in $P'$. 
\end{theorem}

\begin{proof}
Case \eqref{thm:adequate1} is trivial. 
Case \eqref{thm:adequate3} is a consequence of Case \eqref{thm:adequate4}. 
Case \eqref{thm:adequate6} is a consequence of Case \eqref{thm:adequate5}, since every projection is quasi-invertible. 
Cases \eqref{thm:adequate5} and \eqref{thm:adequate7} are a consequence of Case \eqref{thm:adequate2}, for then we already know that $f(P)$ is $\Z$-sup-preserving by Proposition~\ref{prop:baranga}. 
Regarding Cases \eqref{thm:adequate2} and \eqref{thm:adequate4}, we also know that $f(P)$ is $\Z$-sup-preserving by hypothesis or by Proposition~\ref{prop:baranga}, so there it remains to show that $f(P)$ has the $\Z$-refinement property. 
So suppose that $f(x) \leqslant \bigvee_{P'} Z'$ for some $\Z$-subset $Z'$ of $P'$ with sup. 

Case \eqref{thm:adequate2}: 
Since $f$ has the $\Z$-refinement property, there exists some $\Z$-subset $Z$ of $P$ with sup such that $f(Z) \subseteq { \downarrow_{P'}\!\! Z' }$ and $x \leqslant \bigvee_P Z$. 
Now the corestriction $f^{\circ} : P \to f(P)$ is $\Z$-sup-preserving by Lemma~\ref{lem:image}, so $f(x) \leqslant f(\bigvee_P Z) = \bigvee_{f(P)} f(Z)$.
The subset $f(Z)$, as the image of the $\Z$-subset $Z$ of $P$ by $f^{\circ}$, is a $\Z$-subset of $f(P)$, included in $\downarrow_{P'}\!\! Z'$. 
This proves that $f(P)$ has the $\Z$-refinement property in $P'$. 

Case \eqref{thm:adequate4}: 
We have 
$
f(x) = f^{\circ}(g(f(x))) \leqslant \bigvee_{f(P)} f(g(Z'))
$, 
since both $f$ and $g$ are $\Z$-sup-preserving and $f^{\circ}$ is $\Z$-sup-preserving by Lemma~\ref{lem:image}. 
%Moreover, $\bigvee_{P'} f(g(Z')) = f(\bigvee_P g(Z')) \in f(P)$ and so $f(x) \leqslant \bigvee_{f(P)} f(g(Z'))$. 
The inequality $f(g(x')) \leqslant x'$ for all $x' \in P'$ implies that ${ f(g(Z')) } \subseteq { \downarrow_{P'}\!\! Z' }$. 
The subset $f(g(Z'))$, as the image of the $\Z$-subset $g(Z')$ of $P$ by $f^{\circ}$, is a $\Z$-subset of $f(P)$, included in $\downarrow_{P'}\!\! Z'$. 
This proves that $f(P)$ has the $\Z$-refinement property in $P'$. 
\end{proof}

%\begin{remark}
%Note that $f(P)$ is $\Z$-projected in the following cases: 
%\begin{itemize}
%  \item in Case~\eqref{thm:adequate1} (trivially);
%  \item in Case~\eqref{thm:adequate3} if $g$ preserves the $\Z$-below relation, for then we can write $f(P) = p(P')$ where $p$ is the $\Z$-below preserving projection on $P'$ given by $p = f \circ g$;
%  \item in Case~\eqref{thm:adequate5};
%  \item in Case~\eqref{thm:adequate6}.
%\end{itemize}
%\end{remark}

%This is the case since $f(P)$ can be written as $p(P')$, where $p$ is defined by $p = f \circ g : P' \to P'$, and $p$ is a $\Z$-below preserving projection by Proposition~\ref{prop:left}. 

\subsection{Direct images of $\Z$-continuous posets}\label{sec:left2}

In this subsection, we examine conditions under which the direct image $f(P)$ of a $\Z$-continuous poset $P$ by an order-preserving map $f : P \to P'$ is itself $\Z$-continuous, without necessarily being $\Z$-adequate in $P'$ or a $\Z$-subposet of $P'$. 
%We prove new results of this kind and recall existing results by Venugopalan and Baranga.

\begin{definition}
A map $f : P \to P'$ is \textit{weakly $\Z$-adequate} if it is $\Z$-sup-preserving and such that, for all $x' \in f(P)$, there is some $x \in P$ such that $f(x) = x'$ and $f(\twoheaddownarrow^{\Z}_{P} x) \subseteq { \twoheaddownarrow^{\Z}_{P'} x' }$. 
\end{definition}

Using Proposition~\ref{prop:refinement}, it is obvious that every $\Z$-adequate map is weakly $\Z$-adequate. 
We first show the usefulness of this new concept with the following theorem, and then give application cases with Proposition~\ref{prop:wadequ} below.

\begin{theorem}\label{thm:wadequ}
Let $P$ be a $\Z$-continuous poset with weak $\Z$-basis $B$. 
Let $f : P \to P'$ be an order-preserving map with a weakly $\Z$-adequate corestriction $f^{\circ}$. 
Then $f(P)$ is a $\Z$-continuous poset with weak $\Z$-basis $f(B)$. 
Moreover, if $P$ is $\Z$-interpolating and there is some $\Z$-sup-preserving map $g : f(P) \to P$ such that 
\begin{align}\label{eq:swza}
{ \downarrow_{f(P)}\!\! f(\twoheaddownarrow^{\Z}_{P} g(x')) } = { \twoheaddownarrow^{\Z}_{f(P)} x' },
\end{align}
for all $x' \in f(P)$, then $f(P)$ is $\Z$-interpolating. 
\end{theorem}

\begin{proof}
Let $x' \in f(P)$. 
Since $f^{\circ}$ is weakly $\Z$-adequate, there is some $x \in P$ with $f(x) = x'$ and $f(\twoheaddownarrow^{\Z}_{P} x) \subseteq { \twoheaddownarrow^{\Z}_{f(P)} x' }$. 
Since $P$ is $\Z$-continuous with weak $\Z$-basis $B$, there is a $\Z$-subset $Z$ of $P$ included in $\downarrow_{P}\!\!(\twoheaddownarrow^{\Z}_{P} x \cap B)$ with sup $x$. 
The corestriction $f^{\circ}$ of $f$ is $\Z$-sup-preserving by hypothesis, so $x' = \bigvee_{f(P)} f(Z)$, and $f(Z)$ is a $\Z$-subset of $f(P)$ as the direct image of $Z$ by $f^{\circ}$. 
Moreover, 
\begin{align*}
{ f(Z) } &\subseteq { \downarrow_{P'}\!\! f(\twoheaddownarrow^{\Z}_{P} x \cap B) } \subseteq { \downarrow_{P'}\!\! \left(f(\twoheaddownarrow^{\Z}_{P} x) \cap f(B) \right) } \\
&\subseteq { \downarrow_{P'}\!\! \left(\twoheaddownarrow^{\Z}_{f(P)} x' \cap f(B) \right) }. 
\end{align*}
%Since $f(P)$ is $\Z$-sup-preserving in $P'$, we have ${ \twoheaddownarrow^{\Z}_{P'} x' } \subseteq { \twoheaddownarrow^{\Z}_{f(P)} x' }$, hence $f(Z)$ is included in $\downarrow_{P'}\!\! (\twoheaddownarrow^{\Z}_{f(P)} x' \cap f(B))$. 
This proves that $f(P)$ is a $\Z$-continuous poset with weak $\Z$-basis $f(B)$. 

Suppose moreover that $P$ is $\Z$-interpolating and $f$ satisfies Condition~\eqref{eq:swza}, and let $x', z' \in f(P)$ such that $x' \ll^{\Z}_{f(P)} z'$. 
Then, by \eqref{eq:swza}, there is some $y \in P$ such that $x' \leqslant f(y)$ and $y \ll^{\Z}_{P} g(z')$. 
Since $P$ is $\Z$-interpolating, there is some $w \in P$ such that $y \ll^{\Z}_{P} w \ll^{\Z}_{P} g(z')$. 
Now, we have proved that $f(P)$ is $\Z$-continuous, so that $z' = \bigvee_{f(P)} Z'$, for some $\Z$-subset $Z'$ of $f(P)$ included in $\twoheaddownarrow^{\Z}_{f(P)} z'$. 
%Since $f(P)$ is $\Z$-sup-preserving, we have $z' = \bigvee_{P'} Z'$. 
Using that $g : f(P) \to P$ is $\Z$-sup-preserving, we deduce that $g(z') = \bigvee_{P} g(Z')$. 
With $w \ll^{\Z}_{P} g(z')$, this implies that $w \leqslant g(y')$, for some $y' \in Z'$. 
On the one hand, $y \ll^{\Z}_{P} g(y')$, which yields $f(y) \ll^{\Z}_{f(P)} y'$ by \eqref{eq:swza}, so that $x'  \ll^{\Z}_{f(P)} y'$. 
On the other hand, $y' \in { Z' } \subseteq { \twoheaddownarrow^{\Z}_{f(P)} z' }$, so that $y' \ll^{\Z}_{f(P)} z'$. 
This shows that $x' \ll^{\Z}_{f(P)} y' \ll^{\Z}_{f(P)} z'$. 
So $f(P)$ is $\Z$-interpolating, as required. 
\end{proof}

%\begin{remark}
%Be aware that, under the hypotheses of the previous theorem, $f(P)$ does not need to be $\Z$-adequate in $P'$ nor to be a $\Z$-subposet of $P'$. 
%\end{remark}

%A \textit{closure operator} on a poset $P$ is a map $c : P \to P$ such that, for all $x, y \in P$: 
%\begin{itemize}
%	\item $c(x) \leqslant c(y)$ if $x \leqslant y$;
%    \item $x \leqslant c(x)$;
%    \item $c(c(x)) = c(x)$. 
%\end{itemize}
%For instance, if $(f, g)$ is a Galois connection, then $g \circ f$ is a closure operator. 
%The reader may refer e.g.\ to Ern\'e \cite{Erne09} for links between closure operators and Galois connections. 

\begin{proposition}\label{prop:wadequ}
Let $f : P \to P'$ be an order-preserving map. 
Assume that one of the following conditions holds:
\begin{enumerate}
  \item\label{prop:wadequ1bis} $f$ is $\Z$-adequate and surjective;
  \item\label{prop:wadequ1} $f$ is $\Z$-adequate and $f(P)$ is $\Z$-sup-preserving in $P'$;
  \item\label{prop:wadequ2} $f$ is $\Z$-sup-preserving, surjective, and the right adjoint of a Galois connection;
  \item\label{prop:wadequ4} $P$ is cond.\ $\Z$-complete and $f$ is a $\Z$-sup-preserving projection;
  \item\label{prop:wadequ3} $f$ is $\Z$-sup-preserving, $f(P)$ is $\Z$-sup-preserving in $P'$, and there is a $\Z$-sup-preserving map $g : P' \to P$ such that $f \circ g \circ f = f$;
%  \item\label{prop:wadequ5} $P$ is cond.\ $\Z$-complete and $f$ is a $\Z$-sup-preserving closure operator (with $P = P'	$);
  \item\label{prop:wadequ6} $f$ is a $\Z$-sup-preserving order-embedding.
\end{enumerate}
Then the corestriction $f^{\circ}$ of $f$ is weakly $\Z$-adequate. 
Moreover, if $P$ is $\Z$-continuous, then $f^{\circ}$ satisfies Condition~\eqref{eq:swza} of Theorem~\ref{thm:wadequ} for Cases \eqref{prop:wadequ2}, \eqref{prop:wadequ4}, \eqref{prop:wadequ3}, and \eqref{prop:wadequ6}. 
\end{proposition}

\begin{proof}
Case \eqref{prop:wadequ1bis} is a consequence of Case \eqref{prop:wadequ1}. 

Case \eqref{prop:wadequ1}: 
By Lemma~\ref{lem:image}, $f^{\circ}$ is $\Z$-sup-preserving. 
Using the hypothesis that $f(P)$ is $\Z$-sup-preserving in $P'$, it is easy to show that $f^{\circ}$ has the $\Z$-refinement property. 
So $f^{\circ}$ is $\Z$-adequate, hence weakly $\Z$-adequate. 

Case \eqref{prop:wadequ2} is a consequence of Case \eqref{prop:wadequ3}: take for $g$ the left adjoint of $f$.  

Case \eqref{prop:wadequ4} is a consequence of Case \eqref{prop:wadequ3}: take $g := f$, and note that $f(P)$ is $\Z$-sup-preserving in $P'$ by Proposition~\ref{prop:baranga}\eqref{prop:baranga4}.  

Case \eqref{prop:wadequ6}: 
Here, $f^{\circ}$ satisfies the conditions of Case \eqref{prop:wadequ2}, hence is weakly $\Z$-adequate. 

Case \eqref{prop:wadequ3}: 
Let $x' \in f(P)$. 
We show that ${ f(\twoheaddownarrow^{\Z}_{P} x) } \subseteq { \twoheaddownarrow^{\Z}_{f(P)} x' }$ with $x := g(x')$. 
So let $y \ll^{\Z}_{P} x$, and let $x' \leqslant \bigvee_{f(P)} Z'$, for some $\Z$-subset $Z'$ of $f(P)$. 
Since $f(P)$ is $\Z$-sup-preserving in $P'$, $\bigvee_{f(P)} Z' = \bigvee_{P'} Z'$. 
Moreover, $g$ is $\Z$-sup-preserving, so $x = g(x') \leqslant g(\bigvee_{P'} Z') = \bigvee_{P} g(Z')$. 
Now $g(Z')$ is a $\Z$-subset of $P$ and $y \ll^{\Z}_{P} x$, so $y \leqslant g(z')$, for some $z' \in Z'$. 
Thus, $f(y) \leqslant f(g(z')) = z'$. 
Since $f$ is $\Z$-sup-preserving, $f^{\circ}$ is also $\Z$-sup-preserving by Lemma~\ref{lem:image}. 
So $f^{\circ}$ is weakly $\Z$-adequate.

Suppose moreover that $P$ is $\Z$-continuous. 
Then $x = \bigvee_{P} \twoheaddownarrow^{\Z} x$, so $x' = f(x) = f(\bigvee_{P} \twoheaddownarrow^{\Z} x) = \bigvee_{f(P)} f(\twoheaddownarrow^{\Z} x)$. 
This yields ${ \twoheaddownarrow^{\Z}_{f(P)} x' } \subseteq { \downarrow_{f(P)}\!\! f(\twoheaddownarrow^{\Z}_{P} x) }$. 
Note also that $x = h(x')$, where the inclusion map $i : f(P) \to P'$ is $\Z$-sup-preserving by hypothesis and $h := g \circ i : f(P) \to P$ is $\Z$-sup-preserving. 
So $f^{\circ}$ satisfies Condition~\eqref{eq:swza} of Theorem~\ref{thm:wadequ}. 
%Let $x' \in P'$, and let $x := e(x')$, where $e : P' \to P$ is the left adjoint of $f$. 
%Let $y \ll^{\Z}_{P} x$, and let us show that $f(y) \ll^{\Z}_{P'} x'$. 
%So suppose that $x' \leqslant \bigvee_{P'} Z'$, for some $\Z$-subset $Z'$ of $P'$. 
%Then $x = e(x') \leqslant e(\bigvee_{P'} Z') = \bigvee_{P} e(Z')$, since $e$ is $\Z$-sup-preserving by Lemma~\ref{lem:galois}. 
%Now $e(Z')$ is a $\Z$-subset of $P$ and $y \ll^{\Z}_{P} x$, so $y \leqslant e(z')$, for some $z' \in Z'$. 
%So $f(y) \leqslant f(e(z')) = z'$, the latter equality coming from the surjectivity of $f$. 
%This shows that $f(\twoheaddownarrow^{\Z}_{P} x) \subseteq { \twoheaddownarrow^{\Z}_{P'} x' }$. 
%Since $f$ is $\Z$-sup-preserving by hypothesis and $f = f^{\circ}$, we conclude that $f^{\circ}$ is  weakly $\Z$-adequate.
\end{proof}

\begin{remark}
Case \eqref{prop:wadequ2} of the previous proposition, in combination with Theorem~\ref{thm:wadequ}, extends Bandelt and Ern\'e \cite[Theorem]{Bandelt83} and Venugopalan \cite[Theorem~2.17]{Venugopalan86}; it is also close to \cite[Theorem~3.10]{Venugopalan86}, although this latter result was formulated with $\Z$-bases instead of weak $\Z$-bases. 
Similarly, Case \eqref{prop:wadequ4} combined with Theorem~\ref{thm:wadequ} extends Venugopalan \cite[Theorem 2.15]{Venugopalan86} and Baranga \cite[Theorem~2.7]{Baranga96}. 
\end{remark}

\begin{example}\cite{Poncet12a}
Let $(S, \cdot)$ be an inverse semigroup, that is a semigroup such that, for every $x \in S$, there exists a unique $x^* \in S$ with $x x^* x = x$ and $x^* x x^* = x^*$. 
Let $E(S)$ be the commutative subsemigroup made of its idempotent elements. 
Recall that $S$ is naturally equipped with a partial order $\leqslant$ defined by $x \leqslant y$ if $x = y e$ for some $e \in E(S)$. 
We suppose here that $S$ is \textit{mirror}, in the sense that every directed subset of $E(S)$ with a sup in $E(S)$ also has a sup in $S$ \cite[Definition~3.1]{Poncet12a}. In this case, both suprema coincide and belong to $E(S)$, so this amounts to saying that $E(S)$ is a $\D$-sup-preserving subset of $S$. 
The map $\sigma : S \to S, x \mapsto x^* x$ is a $\D$-sup-preserving projection by \cite[Lemma~3.5]{Poncet12a}, whose image is $E(S)$. 
Consequently, $\sigma$ is weakly $\D$-adequate by Proposition~\ref{prop:wadequ}\eqref{prop:wadequ3}.  Hence, by Theorem~\ref{thm:wadequ}, $E(S)$ is a $\D$-continuous poset if $S$ is itself a $\D$-continuous poset, as already stated by \cite[Lemma~5.2]{Poncet12a}. 
\end{example}

The following result extends Venugopalan \cite[Theorem 2.16]{Venugopalan86} and \cite[Theorem 2.17]{Venugopalan86}. 

\begin{corollary}\label{coro:venugo}
Let $P'$ be a %cond.\ $\Z$-complete 
(strongly) $\Z$-continuous poset with weak $\Z$-basis $B'$, and let $f : P \to P'$ be the right adjoint of a Galois connection $(e, f)$. 
Moreover, assume that $f$ is either $\Z$-sup-preserving or a projection (with $P = P'$ in this case). 
Then $f(P)$ is a (strongly) $\Z$-continuous poset with weak $\Z$-basis $f(e(B'))$. 
\end{corollary}

\begin{proof}
By Lemma~\ref{lem:galois}, $e$ is $\Z$-sup-preserving. 
And by Proposition~\ref{prop:baranga}\eqref{prop:baranga1}, $e(P')$ is $\Z$-sup-preserving in $P$. If $f$ is $\Z$-sup-preserving, then $e^{\circ}$ is weakly $\Z$-adequate by Proposition~\ref{prop:wadequ}\eqref{prop:wadequ3}. 
If $f$ is a projection, then $e$ is also a projection, so $e^{\circ}$ is weakly $\Z$-adequate by Proposition~\ref{prop:wadequ}\eqref{prop:wadequ3} (note that we do not use here the weakest assertion of Proposition~\ref{prop:wadequ}\eqref{prop:wadequ4}, since we do not assume $P'$ to be conditionally $\Z$-complete). 
In both cases, we can apply Theorem~\ref{thm:wadequ} to $e$, and we obtain that $e(P')$ is a (strongly) $\Z$-continuous poset with weak $\Z$-basis $e(B')$. 
To conclude, it suffices to note that $e(P')$ and $f(P)$ are order-isomorphic. 
\end{proof}

The following result slightly extends Venugopalan \cite[Theorem 2.18]{Venugopalan86}.
It is also close to \cite[Theorem~3.11]{Venugopalan86}, although the latter was formulated with $\Z$-bases instead of weak $\Z$-bases. 

\begin{corollary}\label{coro:venugo2}
Let $P$ be a %cond.\ $\Z$-complete 
(strongly) $\Z$-continuous poset with weak $\Z$-basis $B$, and let $A$ be a $\Z$-sup-preserving subset of $P$. 
Moreover, assume that the inclusion map $i : A \to P$ is the right adjoint of a Galois connection $(e, i)$. 
Then $A$ is a (strongly) $\Z$-continuous poset with weak $\Z$-basis $e(B)$. 
\end{corollary}

%\begin{corollary}
%Let $P$ be a conditionally $\Z$-complete (strongly) $\Z$-continuous poset with weak $\Z$-basis $B$, and let $c : P \to P$ be a $\Z$-sup-preserving closure operator on $P$. 
%Then $c(P)$ is a (strongly) $\Z$-continuous poset with weak $\Z$-basis $c(B)$. 
%\end{corollary}

%This result is a straightforward consequence of the previous theorem, yet we provide an independent proof. 

%\begin{proof}
%Consider the corestriction $c^{\circ} : P \to c(P)$. 
%This map is surjective. 
%It is also $\Z$-sup-preserving by Lemma~\ref{lem:image}. 
%Moreover, $(c^{\circ}, i)$ is a Galois connection, if we denote by $i$ the inclusion map $c(P) \to P$. 
%By Proposition~\ref{prop:baranga}\eqref{prop:baranga2}, $i$ is a $\Z$-sup-preserving right adjoint. 
%So by Proposition~\ref{prop:left}, $c^{\circ}$ has the $\Z$-refinement property. 
%Thus, Theorem~\ref{thm:imageofzcp} applies, and $c(P)$ is a $\Z$-continuous (or strongly $\Z$-continuous) poset with weak $\Z$-basis $c(B)$. 
%\end{proof}

We are ready to extend Furber's lemma \cite[Lemma~III.1]{Furber19} with the use of the previous corollary. 
A subset $A$ of a poset $P$ is called \textit{$\Z$-inf-preserving} if every $\Z$-subset $Z$ of $A$ with inf in $A$ also has an inf in $P$, and $\bigwedge_P Z = \bigwedge_A Z$. 

\begin{corollary}\label{coro:furber}
Let $P$ be a %cond.\ $\Z$-complete 
(strongly) $\Z$-continuous poset, and let $A$ be a $\Z$-sup-preserving, $\A^*$-inf-preserving subset of $P$ which is also a complete lattice. 
Then $A$ is a (strongly) $\Z$-continuous poset. 
\end{corollary}

\begin{proof}
Take ${ Q } := { \downarrow_{P}\!\! A }$. 
As a lower subset of $P$, $Q$ is $\Z$-adequate in $P$ by Theorem~\ref{thm:lower}, hence is a (strongly) $\Z$-continuous poset by Theorem~\ref{thm:main}. 
%Moreover, it is easily seen that $Q$ is conditionally $\Z$-complete. 
Let us show that $A$ is $\Z$-sup-preserving in $Q$ and that the inclusion map $i : A \to Q$ is the right adjoint of a Galois connection. 
Then, by Corollary~\ref{coro:venugo2}, we will conclude that $A$ is a (strongly) $\Z$-continuous poset. 
So let $Z$ be a $\Z$-subset of $A$ with sup $a_0$ in $A$. 
Since $A$ is $\Z$-sup-preserving in $P$, we have $a_0 = \bigvee_P Z$. 
We have $Z \subseteq Q$ and $\bigvee_P Z = a_0 \in Q$, so $\bigvee_P Z = \bigvee_Q Z$. 
This shows that $a_0 = \bigvee_Q Z$. 
So $A$ is $\Z$-sup-preserving in $Q$. 
Now let $e_A : Q \to A$, $x \mapsto \bigwedge_A (\uparrow_Q\!\! x \cap A)$ (recall that $A$ is a complete lattice by hypothesis), and let $x \in Q$, $a \in A$. 
If $x \leqslant i(a) = a$, then $a \in \uparrow_Q\!\! x \cap A$, so $e_A(x) \leqslant a$. 
Conversely, suppose that $e_A(x) \leqslant a$. 
Since $A$ is $\A^*$-inf-preserving in $P$ and $\uparrow_Q\!\! x \cap A$ is nonempty by definition of $Q$, we have $e_A(x) = \bigwedge_P (\uparrow_Q\!\! x \cap A)$, which implies that $x \leqslant e_A(x)$. 
Thus, $x \leqslant a = i(a)$. 
So $(e_A, i)$ is a Galois connection, as required.
\end{proof}

\subsection{Direct images of subsets}\label{sec:left3}

In this last subsection, we consider subsets of the form $f(A)$, for some subset $A$ of a poset $P$ and some order-preserving map $f : P \to P'$. 

%\begin{lemma}
%Let $A$ be a $\Z$-adequate subset of a poset $P$, and let $f : P \to P'$ be a $\Z$-adequate map. 
%Then $f(A)$ has the $\Z$-refinement property in $P'$. 
%\end{lemma}

%\begin{proof}
%Let $i : A \to P$ be the inclusion map. 
%Since $A$ is $\Z$-sup-preserving, $i$ is $\Z$-sup-preserving, so $f \circ i$ is $\Z$-sup-preserving. 
%Moreover, both $i$ and $f$ have the $\Z$-refinement property, so $f \circ i$ has the $\Z$-refinement property. 
%So $f \circ i$ is $\Z$-adequate. 
%From Proposition~\ref{prop:refin}\eqref{prop:refin2} we deduce that $f \circ i(A) = f(A)$ has the $\Z$-refinement property. 
%\end{proof}

\begin{theorem}\label{thm:image}
Let $A$ be a $\Z$-adequate subset of a poset $P$, and let $f : P \to P'$ be the left adjoint of a Galois connection $(f, g)$ where $g$ is $\Z$-sup-preserving. 
Assume that one of the following conditions holds:
\begin{enumerate}
  \item\label{thm:image1} $f(g(P')) \subseteq { f(A) }$; 
  \item\label{thm:image2} $g(f(a)) = a$ for all $a \in A$ (this holds if $f$ is injective); 
  \item\label{thm:image3} $A$ is $\A$-sup-regular and $P$ is conditionally $\A$-complete; 
  \item\label{thm:image4} $A$ is $\Z$-sup-regular with $g(f(A)) \subseteq { A }$ and $P$ is conditionally $\Z$-complete.  
\end{enumerate}
Then $f(A)$ is $\Z$-adequate in $P'$. 
\end{theorem}

\begin{proof}
Let $i : A \to P$ be the inclusion map, which is $\Z$-adequate. 
The map $f$ being $\Z$-adequate by Proposition~\ref{prop:left}, the composition $f \circ i$ is also $\Z$-adequate. 
Hence, considering Theorem~\ref{thm:adequate}\eqref{thm:adequate2}, it suffices to show that $f(A)$ is $\Z$-sup-preserving in $P'$ in order to prove the result. 
So let $Z'$ be a $\Z$-subset of $f(A)$ with a sup in $f(A)$. 
We write this sup by $f(a_0)$ with $a_0 \in A$. %, and since $f(a_0) = f(g(f(a_0)))$ we can suppose that $a_0 = g(f(a_0))$ without loss of generality. This implies that $a_0$ is an upper bound of $g(Z')$. DERNIERES PHRASES FAUSSES CAR $g(f(a_0))$ N'EST PAS DANS A EN GENERAL
Let $u'$ be an upper bound of $Z'$ in $P'$. 
Then $z' \leqslant u'$ for all $z' \in Z'$. 
We want to show that $f(a_0) \leqslant u'$, i.e.\ $a_0 \leqslant g(u')$. 

Case \eqref{thm:image1}: 
We have $z' = f(g(z')) \leqslant f(g(u')) \leqslant u'$ for all $z' \in Z'$. 
Since $f(g(P')) \subseteq { f(A) }$, there exists some $a \in A$ such that $f(g(u')) = f(a)$. 
So $f(a)$ is an upper bound of $Z'$, hence $f(a_0) \leqslant f(a)$ by definition of $f(a_0)$. 
This shows that $f(a_0) \leqslant u'$, so $f(a_0)$ is the least upper bound of $Z'$ in $P'$, i.e.\ its sup in $P'$. 

Case \eqref{thm:image2}:
We have $g(Z') \subseteq { g(f(A)) } \subseteq { A }$. 
Let us show that $a_0$ is the sup of $g(Z')$ in $A$. 
First, $a_0$ is an upper bound of $g(Z')$, since $g(z') \leqslant g(f(a_0)) = a_0$. 
Let $a$ be another upper bound of $g(Z')$ in $A$. 
Then $f(a)$ is an upper bound of $Z'$, so $f(a_0) \leqslant f(a)$. 
This implies $a_0 \leqslant g(f(a)) = a$. 
So $a_0$ is the sup of $g(Z')$ in $A$. 
Since $A$ is $\Z$-sup-preserving, $a_0$ is also the sup of $g(Z')$ in $P$. 
Notice that $g(u')$ is an upper bound of $g(Z')$, so $a_0 \leqslant g(u')$, as required.

Case \eqref{thm:image3}:
If $z' \in Z'$, there is some $a(z') \in A$ such that $z' = f(a(z'))$. 
Since $P$ is conditionally $\A$-complete, we can define $x = \bigvee_P a(Z')$, and since $A$ is $\A$-sup-regular we have $x \in A$. 
Now $f(x) = \bigvee_{P'} f(a(Z')) = \bigvee_{P'} Z'$. 
This yields $\bigvee_{P'} Z' \in f(A)$, so $f(a_0) = \bigvee_{f(A)} Z' = \bigvee_{P'} Z' \leqslant u'$. 

Case \eqref{thm:image4}: 
Take $Z := g(Z')$. This is a $\Z$-subset of $P$. 
Moreover, since $g(f(A)) \subseteq { A }$ by hypothesis, $Z$ is a subset of $A$. 
Now $P$ is conditionally $\Z$-complete, so we can define $z_0 = \bigvee_P Z$. 
By $\Z$-sup-regularity of $A$ we have $z_0 \in A$. 
We deduce that $f(z_0) = \bigvee_{P'} f(Z) = \bigvee_{P'} f(g(Z'))$. 
Since $Z' \subseteq { f(A) }$ we have $f(g(Z')) = Z'$, so $f(z_0) = \bigvee_{P'} Z'$. 
This proves that $\bigvee_{P'} Z' \in f(A)$, so $f(a_0) = \bigvee_{f(A)} Z' = \bigvee_{P'} Z' \leqslant u'$. 
\end{proof}

\begin{theorem}
Let $A'$ be a $\Z$-adequate subset of a poset $P'$, and let $(f, g)$ be a Galois connection such that $g : P' \to P$ is $\Z$-sup-preserving surjective and $f(g(A')) \subseteq { A' }$. 
Then $g(A')$ coincides with $f^{-1}(A')$ and is $\Z$-adequate in $P$. 
Moreover, if $A'$ is order-convex (resp.\ $\Z$-sup-generates $P'$), then $g(A')$ is order-convex (resp.\ $\Z$-sup-generates $P$). 
\end{theorem}

\begin{proof}
Take $A := g(A')$. 
The equality $A = f^{-1}(A')$ is an easy consequence of the hypothesis $A \subseteq { f^{-1}(A') }$ and the fact that $g$ is surjective.
Let $a \in A$, and let $Z$ be a $\Z$-subset of $P$ with sup such that $a \leqslant \bigvee_P Z$. 
Then $f(a) \leqslant f(\bigvee_P Z) = \bigvee_{P'} f(Z)$. 
Since $f(a) \in A'$, $f(Z)$ is a $\Z$-subset of $P'$, and $A'$ has the $\Z$-refinement property in $P'$, we deduce that there exists some $\Z$-subset $Z'$ of $A'$ such that $Z' \subseteq { \downarrow_{P'}\!\! f(Z) }$ and $f(a) \leqslant \bigvee_{A'} Z'$. 
Now $A'$ is $\Z$-sup-preserving in $P'$, so $f(a) \leqslant \bigvee_{P'} Z'$. 
Using that $g$ is $\Z$-sup-preserving, we get $a \leqslant g(\bigvee_{P'} Z') = \bigvee_P g(Z')$. 
Now, $Z' \subseteq { A' }$ implies that $g(Z') \subseteq { g(A') } = { A }$ and $Z' \subseteq { \downarrow_{P'}\!\! f(Z) }$ implies that ${ g(Z') } \subseteq { g(\downarrow_{P'}\!\!  f(Z)) } \subseteq { \downarrow_{P}\!\! g(f(Z)) } = { \downarrow_{P}\!\!  Z }$, where the last equality holds because $f$ is injective. 
Consequently, $g(Z') \subseteq { \downarrow_{P}\!\! Z \cap A }$. 
Moreover, $\bigvee_P g(Z') = g(\bigvee_{P'} Z') = g(\bigvee_{A'} Z') \in g(A') = A$, so $\bigvee_P g(Z') = \bigvee_A g(Z')$, hence $a \leqslant \bigvee_A g(Z')$. 
Furthermore, $g(Z')$ is a $\Z$-subset of $A$ as the image of $Z'$ by the order-preserving map $A' \to A$, $x' \mapsto g(x')$. 
This proves that $A$ has the $\Z$-refinement property. 

Now we show that $A$ is $\Z$-sup-preserving. 
Let $Z$ be a $\Z$-subset of $A$ with sup $a_0$ in $A$. 
First remark that $f(a_0)$ is the sup of $f(Z)$ in $A'$. 
Indeed, if $u' \in A'$ is an upper bound of $f(Z)$, then $g(u') \in g(A') = { A }$ is an upper bound of $Z$ in $A$, so $a_0 \leqslant g(u')$, thus $f(a_0) \leqslant u'$. 
Since $A'$ is $\Z$-sup-preserving, this implies that $f(a_0)$ is the sup of $f(Z)$ in $P'$. 
To conclude, let us show that $a_0$ is the sup of $Z$ in $P$. 
If $u \in P$ is an upper bound of $Z$, then $f(u)$ is an upper bound of $f(Z)$ in $P'$, so $f(a_0) \leqslant f(u)$. Thus, $a_0 \leqslant g(f(u)) = u$, the latter equality coming from the surjectivity of $g$. 

Suppose that $A'$ is order-convex. 
If $a_1, a_2 \in A$ and $a_1 \leqslant x \leqslant a_2$, then $f(a_1) \leqslant f(x) \leqslant f(a_2)$. 
Since $f(a_1), f(a_2) \in A'$ and $A'$ is order-convex, we have $f(x) \in A'$. 
So $x \in f^{-1}(A') = A$. 
This shows that $A$ is order-convex. 

Suppose that $A'$ $\Z$-sup-generates $P'$. 
Let $x \in P$. 
Since $g$ is surjective, there is some $x' \in P'$ with $x = g(x')$. 
Moreover, $A'$ $\Z$-sup-generates $P'$, so there is some $\Z$-subset $Z'$ of $P'$ included in $A'$ such that $x' = \bigvee_{P'} Z'$. 
This implies that $g(Z') \subseteq { g(A') } = { A }$ and $x = g(x') = g(\bigvee_{P'} Z') = \bigvee_P g(Z')$. 
This shows that $A = g(A')$ $\Z$-sup-generates $P$, as required. 
\end{proof}

%%%%%%%%%%%%%%%%%%%%%%
%%%%%%%%%%%%%%%%%%%%%%
%%%%%%%%%%%%%%%%%%%%%%
%%%%%%%%%%%%%%%%%%%%%%
\section{Kernel retracts as $\Z$-subposets}\label{sec:kernels}

%\subsection{Kernels on posets}

A \textit{kernel operator} on a poset $P$ is a map $k : P \to P$ such that, for all $x, y \in P$: 
\begin{itemize}
	\item $k(x) \leqslant k(y)$ if $x \leqslant y$, 
	\item $k(x) \leqslant x$, 
	\item $k(k(x)) = k(x)$. 
\end{itemize}
For instance, if $(f, g)$ is a Galois connection, then $f \circ g$ is a kernel operator. 
Conversely, every kernel operator arises this way. 
To see this, let $k$ be a kernel operator. 
Define $K := \{ k(x) : x \in P \}$, let $i : K \to P$ be the inclusion map, and let $k^{\circ}$ be the corestriction of $k$. 
Then $k = i \circ k^{\circ}$. 
Notice that $K = \{ x \in P : k(x) = x \} = \{ x \in P : k^{\circ}(x) = x \}$. 
The surjective map $k^{\circ}$ is called a \textit{kernel retraction} on $P$ and the subset $K$ of $P$ is called the \textit{kernel retract} of $P$ induced by $k$. 
In the following, we see that $(i, k^{\circ})$ is a Galois connection. 
%A \textit{kernel retraction} on $P$ is a surjective map $k^{\circ} : P \to K$ such that the map $P \ni x \mapsto k^{\circ}(x) \in P$ is a kernel operator on $P$. 
%In this case we have $K = \{ x \in P : k^{\circ}(x) = x \}$, and the subset $K$ of $P$ is called a \textit{kernel retract} of $P$ induced by the kernel retraction $k^{\circ}$. 

\begin{lemma}\label{lem:ker}
Every kernel retraction is the right adjoint of a Galois connection. %; in particular it preserves arbitrary existing infima. 
%Consequently, it is Scott-continuous if and only if it is Lawson-continuous. 
\end{lemma}

\begin{proof}
Let $k^{\circ} : P \to K$ be a kernel retraction. 
Then $y \leqslant k^{\circ}(x)$ if and only if $y \leqslant x$, for all $x \in P$, $y \in K$. 
Thus, $(i, k^{\circ})$ is a Galois connection, where $i : K \to P$ is the inclusion map. 
%This implies that $k^{\circ}$ preserves arbitrary existing infima (see e.g.\ \cite[Theorem~O-3.3]{Gierz03}). 
% and that, if $k^{\circ}$ is Scott-continuous, it is Lawson-continuous (use \cite[Lemma~III-1.2(i)]{Gierz03}). %Also, if $k^{\circ}$ is Lawson-continuous then it is Scott-continuous by a remark made in Paragraph~\ref{subsec:sl}. 
\end{proof}

\begin{lemma}\label{lem:aregker}
Every kernel retract of a poset is $\A$-sup-regular. 
\end{lemma}

\begin{proof}
Let $k^{\circ} : P \to K$ be a kernel retraction associated with the kernel retract $K$ of a poset $P$, and let $i : K \to P$ be the inclusion map. 
Then $(i, k^{\circ})$ is a Galois connection, and $i$ is $\A$-sup-preserving by Lemma~\ref{lem:galois}, which means that $K$ is $\A$-sup-preserving in $P$. 
Let $A$ be a (possibly empty) subset of $K$ with a sup $a_0$ in $P$. 
Then $k^{\circ}(a_0)$ is an upper bound of $A$ and $k^{\circ}(a_0) \leqslant a_0$, so $k^{\circ}(a_0) = a_0$, hence $a_0 \in K$. 
This implies that $a_0$ is the sup of $A$ in $K$. 
This shows that $K$ is $\A$-sup-regular. 
\end{proof}

\begin{remark}
Let $k : P \to P$ be a kernel operator, and $i : k(P) \to P$ be the inclusion map. 
Combining Lemma~\ref{lem:image} with Lemma~\ref{lem:aregker}, and recalling that $k = i \circ k^{\circ}$, we deduce that $k$ is $\Z$-sup-preserving if and only if its corestriction $k^{\circ} : P \to k(P)$ is $\Z$-sup-preserving. 
See also \cite[Proposition~1.13]{Baranga96}. 
\end{remark}

\begin{theorem}\label{thm:thmQ}
Let $K$ be a kernel retract of a poset $P$ induced by a $\Z$-sup-preserving kernel retraction $k^{\circ} : P \to K$. 
Then $K$ is $\A$-sup-regular and $\Z$-adequate in $P$. 
Moreover, if $P$ is (strongly) $\Z$-continuous with weak $\Z$-basis $B$, then $K$ is (strongly) $\Z$-continuous with weak $\Z$-basis $k(B)$. 
\end{theorem}

\begin{proof}
The subset $K$ is $\A$-sup-regular by Lemma~\ref{lem:aregker}, and is $\Z$-adequate by application of Theorem~\ref{thm:adequate}\eqref{thm:adequate3} to the Galois connection $(i, k^{\circ})$ of Lemma~\ref{lem:ker}. 
Now suppose that $P$ is (strongly) $\Z$-continuous with weak $\Z$-basis $B$. 
By Proposition~\ref{prop:zdense}\eqref{prop:zdense6}, $k(B)$ is $\Z$-dense in $K$ rel.\ $(P, B)$. 
So $K$ is (strongly) $\Z$-continuous with weak $\Z$-basis $k(B)$ by Theorem~\ref{thm:main}\eqref{thm:main3}. 
\end{proof}

%\begin{remark}
%While the second part of the previous theorem may seem to be a special case of Venugopalan's theorem (Theorem~\ref{thm:venugo}), we do not need to require that $P$ be conditionally $\Z$-complete here. 
%\end{remark}

%The following result corresponds to Venugopalan \cite[Theorem~2.17]{Venugopalan86}. 

%\begin{corollary}[Venugopalan]
%Let $Q$ be a (strongly) $\Z$-continuous poset and let $g : Q \to P$ be a $\Z$-sup-preserving surjective right map. 
%Then $P$ is a (strongly) $\Z$-continuous poset. 
%\end{corollary}

%%%%%%%%%%%%%%%%%%%%%%
%%%%%%%%%%%%%%%%%%%%%%
%%%%%%%%%%%%%%%%%%%%%%
%%%%%%%%%%%%%%%%%%%%%%
\section{Existence of a largest $\D$-continuous $\D$-subposet}\label{sec:largest}

In this section, we deduce from previous results the theorem announced in the introduction, which is that every cond.\ $\A$-complete, $\D$-interpolating poset contains a largest $\D$-continuous $\D$-subposet. 

Recall that, in a poset $P$, a $\Z$-ideal is a subset $I$ of $P$ such that $I = { \downarrow\!\! Z }$ for some $\Z$-subset $Z$ of $P$. 
We consider the following subset of a poset $P$: 
\[
P^{\Z} := \{ x \in P : \twoheaddownarrow^{\Z} x \mbox{ is a $\Z$-ideal and has a sup } \}. 
\]

\begin{lemma}\label{lem:knabla}
Let $P$ be a $\Z$-interpolating poset. % such that $P^{\Z}$ is nonempty. 
%If $\twoheaddownarrow x$ has a sup for all $x \in P^\nabla$, 
Then the map $k^{\Z}_P : P^{\Z} \to P^{\Z}$ defined by 
\[
k^{\Z}_P(x) = \bigvee \twoheaddownarrow^{\Z} x, 
\]
for all $x \in P^{\Z}$, is a kernel operator on $P^{\Z}$ such that $\twoheaddownarrow^{\Z} x = \twoheaddownarrow^{\Z} k^{\Z}_P(x)$, for all $x \in P^{\Z}$. 
\end{lemma}

\begin{proof}
Let $x, y \in P^{\Z}$ with $x \leqslant y$. 
From the definition of $k^{\Z}_P$ we easily deduce that $k^{\Z}_P(x) \leqslant x$ and $k^{\Z}_P(x) \leqslant k^{\Z}_P(y)$. 
To prove that $k^{\Z}_P(x) \in P^{\Z}$, let us show first that the following equivalence holds for all $v \in P$: 
\begin{equation}\label{eq:EQ}
{ v \ll^{\Z} x } \Leftrightarrow { v \ll^{\Z} k^{\Z}_P(x) }. 
\end{equation}
So let $v \ll^{\Z} x$. 
Since $P$ is $\Z$-interpolating, $v \ll^{\Z} w \ll^{\Z} x$ for some $w \in P$, and the definition of $k^{\Z}_P$ gives $w \leqslant k^{\Z}_P(x)$, hence $v \ll^{\Z} k^{\Z}_P(x)$. 
Conversely, the assertion $v \ll^{\Z} k^{\Z}_P(x)$ trivially implies $v \ll^{\Z} x$ since $k^{\Z}_P(x) \leqslant x$. 
Equivalence~\eqref{eq:EQ} can be rewritten $\twoheaddownarrow^{\Z} x = \twoheaddownarrow^{\Z} k^{\Z}_P(x)$, hence $k^{\Z}_P(x) \in P^{\Z}$. 
We also get $k^{\Z}_P(k^{\Z}_P(x)) = \bigvee \twoheaddownarrow^{\Z} k^{\Z}_P(x) = \bigvee \twoheaddownarrow^{\Z} x = k^{\Z}_P(x)$. 
This proves that $k^{\Z}_P$ is a kernel operator on $P^{\Z}$. 
\end{proof}

\begin{definition}
A poset $P$ is a \textit{$\Z$-Mao--Xu poset} if $\twoheaddownarrow^{\Z} x$ is a $\Z$-ideal with sup for all $x \in P$, in other words if $P = P^{\Z}$. 
\end{definition}

Every $\Z$-continuous poset is obviously $\Z$-Mao--Xu. 
The following example gives some additional insight in the case where $\Z = \D$. 

\begin{example}\label{ex:starstar}
Let $P$ be a conditionally $\A^*$-complete poset. 
Then $x \in P^{\D}$ if and only if $D_x := \twoheaddownarrow^{\D}_P x$ is nonempty, so that $P^{\D} = \twoheaduparrow{}^{\D}_P P$. 
Indeed, if $x \in P^{\D}$, then $D_x$ is nonempty as a directed subset. 
Conversely, assume that $D_x$ is nonempty. 
Then $D_x$ has a sup since it is bounded above by $x$ and $P$ is conditionally $\A^*$-complete. 
Moreover, $D_x$ is directed, since $v, w \in D_x$ implies $v \vee w \in D_x$ (see e.g.\ \cite[Proposition~I-1.2(iii)]{Gierz03}). 
So $x \in P^{\D}$, as required. 
This shows in particular that $P^{\D}$ is an upper subset of $P$, so $P^{\D}$ is $\A^*$-sup-regular by Theorem~\ref{thm:sup2}. 
From this latter property, we easily deduce that $P^{\D}$ is itself conditionally $\A^*$-complete.  

If $P$ is conditionally $\A^*$-complete and $\twoheaddownarrow^{\D} x$ is nonempty for all $x \in P$, then $x \in P^{\D}$ for all $x \in P$, so $P$ is $\D$-Mao--Xu. 
In particular, if $P$ is (conditionally) $\A$-complete, then $P$ is $\D$-Mao--Xu. 
%In particular $P^{\D}$ is a $\D$-Scott-open subset of $P$. 

If $P$ is conditionally $\A^*$-complete and $\D$-interpolating, then $P^{\D} = \twoheaduparrow{}^{\D}_P P$ as noted above, so $P^{\D}$ is a $\D$-Scott-open subset of $P$ by Example~\ref{ex:scottopen}. 
By Theorem~\ref{thm:sup2}, $P^{\D}$ is $\A^*$-sup-regular and $\D$-adequate. 
So, by Theorem~\ref{thm:main}, $P^{\D}$ is a $\D$-subposet of $P$. 
This implies, by Theorem~\ref{thm:main}\eqref{thm:main1}, that $P^{\D}$ is $\D$-interpolating as an upper subset of $P$. 
Moreover, %$P^{\D}$ is $\D$-interpolating and conditionally $\A^*$-complete, and 
\[
(P^{\D})^{\D} = P^{\D}, 
\]
i.e.\ $P^{\D}$ is a $\D$-Mao--Xu poset. 
Indeed, let $x \in P^{\D}$. 
From Lemma~\ref{lem:knabla} we have $\bigvee_P D_x = k^{\D}_P(x) \in P^{\D}$, this latter subset being $\D$-Scott-open. 
Thus, $D_x \cap P^{\D}$ is nonempty. 
Using the fact that $P^{\D}$ is a $\D$-subposet of $P$, we have $D_x \cap P^{\D} = \twoheaddownarrow^{\D}_{P^{\D}} x$. 
As remarked above, $P^{\D}$ is itself conditionally $\A^*$-complete, so the nonemptiness of $\twoheaddownarrow^{\D}_{P^{\D}} x$ is sufficient for $x$ to be in $(P^{\D})^{\D}$.
\end{example}

%Recall from Example~\ref{ex:starstar} that a poset $P$ is $\D$-Mao--Xu in any of the following situations: 
%\begin{itemize}
%    \item $P$ is conditionally $\A^*$-complete and $\twoheaddownarrow^{\D} x$ is nonempty for all $x$; 
%	\item $P$ is conditionally $\A$-complete;
%	\item $P$ is $\A$-complete; 
%	\item $P$ is a $\D$-continuous poset. 
%\end{itemize}

The following result extends \cite[Theorem~3.11]{Mao17}. 

\begin{corollary}
If $P$ is a $\Z$-Mao--Xu poset, then $P$ is strongly $\Z$-continuous if and only if $P$ is $\Z$-interpolating and $\twoheaddownarrow^{\Z} x = \twoheaddownarrow^{\Z} y$ implies $x = y$ for all $x, y \in P$. 
\end{corollary}

\begin{proof}
The `only if' part is clear. 
For the `if' part, let $x \in P$. 
Since $P$ is a $\Z$-Mao--Xu poset, then $P = P^{\Z}$. 
By the previous lemma, $\twoheaddownarrow^{\Z} x = \twoheaddownarrow^{\Z} k^{\Z}_P(x)$. 
By hypothesis we get $x = k^{\Z}_P(x)$, i.e.\ $x$ is the sup of $\twoheaddownarrow^{\Z} x$. 
We already know that $\twoheaddownarrow^{\Z} x$ is a $\Z$-ideal. 
So $P$ is (strongly) $\Z$-continuous. 
\end{proof}

%Note that, if $P$ is conditionally-complete, then $P$ is Mao--Xu if and only if $\twoheaddownarrow x$ is nonempty for all $x$. 
%$\twoheaddownarrow x$ is directed if and only if $\twoheaddownarrow x$ is nonempty. %, and this is always the case if $P$ is continuous. %A \textit{domain} is a continuous poset in which every directed subset has a sup. 

\begin{lemma}
Let $P$ be a $\Z$-interpolating poset such that $P^{\Z}$ is a $\Z$-sup-preserving subset. 
Then $k^{\Z}_P(P^{\Z})$ is a kernel retract of $P^{\Z}$ induced by a $\Z$-sup-preserving kernel retraction. 
\end{lemma}

\begin{proof}
Let $k^{\Z}_P : P^{\Z} \to P^{\Z}$ be the kernel operator on $P^{\Z}$ of Lemma~\ref{lem:knabla}. 
Take $K := k^{\Z}_P(P^{\Z}) = \{ x \in P^{\Z} : k^{\Z}_P(x) = x\}$, and let $k_P^{\circ}$ be the corestriction of $k^{\Z}_P$. 
Then $k_P^{\circ}$ is a kernel retraction, so $K$ is a kernel retract of $P^{\Z}$. 
Let us show that $k_P^{\circ}$ is $\Z$-sup-preserving. 
This amounts to show that, if $Z$ is a $\Z$-subset of $P^{\Z}$ with sup $z_0$ in $P^{\Z}$, then $k^{\Z}_P(z_0) = \bigvee_K k^{\Z}_P(Z)$. 
Obviously, $k^{\Z}_P(z_0)$ is an upper bound of $k^{\Z}_P(Z)$ in $K$.
Let $u$ be another upper bound of $k^{\Z}_P(Z)$ in $K$. 
Since $P^{\Z}$ is supposed to be $\Z$-sup-preserving, $z_0$ is also the sup of $Z$ in $P$. 
Now if $x \ll^{\Z}_P z_0 = \bigvee_P Z$, there is some $z \in Z$ such that $x \ll^{\Z}_P z$ (using the definition of the $\Z$-below relation and the $\Z$-interpolation property). 
This implies that $x \ll^{\Z}_P k^{\Z}_P(z)$ by the equivalence \eqref{eq:EQ} in the proof of Lemma~\ref{lem:knabla}, hence $x \leqslant u$. 
Using the definition of $k^{\Z}_P$ we conclude that $k^{\Z}_P(z_0) \leqslant u$. 
Thus, $k^{\Z}_P(z_0)$ is the least upper bound of $k^{\Z}_P(Z)$ in $K$, i.e.\ $k^{\Z}_P(z_0) = \bigvee_K k^{\Z}_P(Z)$, as required. 
\end{proof}

For the following theorem we need to specialize to the case $\Z = \D$. 

\begin{theorem}
%Let $P$ be a conditionally-complete interpolating poset such that $P^\nabla \neq \emptyset$. 
Let $P$ be a $\D$-interpolating poset such that $P^{\D}$ is a $\D$-Scott-open subset. 
%If $\twoheaddownarrow x$ has a sup for all $x \in P^\nabla$, 
Then $k^{\D}_P(P^{\D})$ is a kernel retract of $P^{\D}$ induced by a $\D$-sup-preserving kernel retraction. 
Moreover, this is the largest $\D$-continuous $\D$-subposet of $P$. 
% admits a maximal continuous kernel retract. 
%Then $\im k$ is a kernel retract of $P$, and a conditionally-complete continuous poset. 
\end{theorem}

\begin{proof}
%From Theorem~\ref{thm:thmQ} we deduce 
In this proof we write $U$ for the $\D$-Scott-open subset $P^{\D}$, and $k$ for the map $k^{\D}_P$. 
Take $K := k(U) = \{ x \in U : k(x) = x\}$. 
We already know that $U$ is $\D$-sup-preserving by Theorem~\ref{thm:sup2}, so $k$ is a $\D$-sup-preserving kernel retraction and $K$ is a kernel retract of $U$ by the previous lemma. 
By Theorem~\ref{thm:thmQ}, $K$ is $\D$-adequate in $P$, hence is a $\D$-subposet of $P$ by Theorem~\ref{thm:main}. 
Now let us show that $K$ is a $\D$-continuous poset, so let $x \in K$. 
We need to show that $\twoheaddownarrow^{\D}_K x = \twoheaddownarrow^{\D}_P x \cap K$ contains a directed subset that admits $x$ as sup in $K$. 

We have $\bigvee_P \twoheaddownarrow^{\D}_P x = k(x) = x \in U$; since $U$ is $\D$-Scott-open we deduce that $\twoheaddownarrow^{\D}_P x \cap U$ is nonempty. 
Applying Lemma~\ref{lem:sup} to the directed subset $\twoheaddownarrow^{\D}_P x$, we get that $\twoheaddownarrow^{\D}_P x \cap U$ is directed and $\bigvee_P \twoheaddownarrow^{\D}_P x = \bigvee_P (\twoheaddownarrow^{\D}_P x \cap U)$. 
So $x = k(x) = k(k(x)) = k(\bigvee_P \twoheaddownarrow^{\D}_P x) = k(\bigvee_P (\twoheaddownarrow^{\D}_P x \cap U))$. 
Since the corestriction $k^{\circ}$ is $\D$-sup-preserving by Lemma~\ref{lem:image}, we obtain $x = \bigvee_K k(\twoheaddownarrow^{\D}_P x \cap U)$. 
Thus, $x$ is the sup in $K$ of the (non-empty) directed subset $k(\twoheaddownarrow^{\D}_P x \cap U)$, which is itself included in $\twoheaddownarrow^{\D}_P x \cap K = \twoheaddownarrow^{\D}_K x$. 
This proves that $K$ is a $\D$-continuous poset. 

To conclude the proof we show that $K$ is the largest $\D$-continuous $\D$-subposet of $P$. 
So let $L$ be a $\D$-continuous $\D$-subposet of $P$, and let us show that $L \subseteq { K }$. 
Let $y \in L$. 
Then $\twoheaddownarrow^{\D}_L y = \twoheaddownarrow^{\D}_P y \cap L$, and this subset is directed. 
Since $L$ is $\D$-continuous, $y = \bigvee_L \twoheaddownarrow^{\D}_L y$. 
Moreover, every $\D$-subposet is $\D$-sup-preserving by definition, so $y = \bigvee_P \twoheaddownarrow^{\D}_P y \cap L$. 
This implies that $\twoheaddownarrow^{\D}_P y = \downarrow\!\! (\twoheaddownarrow^{\D}_P y \cap L)$, hence $y \in U$, and $y = \bigvee_P \twoheaddownarrow^{\D}_P y = k(y)$, i.e.\ $y \in K$. 
This proves that $L \subseteq { K }$. 
\end{proof}

The following result due to Mao and Xu \cite[Theorem~3.8]{Mao17} now comes as a corollary to the previous theorem. 

\begin{corollary}[Mao--Xu]
Let $P$ be a $\D$-interpolating, $\D$-Mao--Xu, $\D$-complete poset. 
Then $k^{\D}_P(P)$ is a $\D$-continuous $\D$-complete poset. 
\end{corollary}

\begin{proof}
Since $P$ is $\D$-Mao--Xu, we have $P^{\D} = P$. 
In particular, $P^{\D}$ is $\D$-Scott-open in $P$. 
Hence, the conditions of the previous theorem are fulfilled, so $k^{\D}_P(P)$ is $\D$-continuous. 
By Lemma~\ref{lem:aregker}, $k^{\D}_P(P)$ is $\A$-sup-regular. 
Since $P$ is supposed to be $\D$-complete, $k^{\D}_P(P)$ is thus also $\D$-complete.  
\end{proof}

\begin{example}[Continuation of Example~\ref{ex:one} and Figure~\ref{fig:Hasse}]
Since $P$ is a complete lattice, we have $P = P^{\D}$. 
Moreover, $k^{\D}_P(\omega) = 0$ and $k^{\D}_P(x) = x$ for every $x \neq \omega$, so that $k^{\D}_P(P) = P \setminus \{ \omega \}$. 
\end{example}

\begin{example}[Continuation of Example~\ref{ex:two} and Figure~\ref{fig:Hasse2}]
Again $P$ is a complete lattice, so $P = P^{\D}$. 
Now in this poset the assertion $x \ll^{\D} y$ never holds, except if $x = 0$. 
Thus, $k^{\D}_P(P) = \{ 0 \}$. 
Note that the poset $P \setminus \{ \omega \}$ is a $\D$-continuous poset (every element is $\D$-compact), but is not a $\D$-subposet of $P$. 
\end{example}

%\begin{proposition}[Bandelt--Ern\'e]
%Let $(f, g)$ be a Galois connection between posets $P$ and $Q$.  
%If $Q$ is (conditionally) $\Z$-complete, then $P$ is (conditionally) $\Z$-complete. 
%Moreover, 
%\[
%\bigvee_P Z = g(\bigvee_Q f(Z)), 
%\]
%for every (upper-bounded) $\Z$-subset $Z$ of $P$. 
%\end{proposition}

%\begin{proof}
%We recall the proof of \cite{Bandelt83} for the sake of completeness. 
%It is clear that $g(\bigvee_Q f(Z))$ is an upper bound of $Z$ in $P$. 
%Let $u$ be another upper bound. 
%Then $f(u)$ is an upper bound of $f(Z)$ in $Q$, so $\bigvee_Q f(Z) \leqslant f(u)$. 
%This shows that $g(\bigvee_Q f(Z)) \leqslant g(f(u)) = u$, so $g(\bigvee_Q f(Z))$ is the sup of $Z$ in $P$. 
%\end{proof}

\begin{acknowledgements}
I would like to gratefully thank Marianne Akian for her useful remarks and advice on a preliminary version of the manuscript. 
I am also indebted in two anonymous referees who made a comprehensive and detailed review and provided me with a great number of valuable suggestions to improve readability, correctness, and consistency with the existing literature. 
\end{acknowledgements}

\bibliographystyle{plain}

\def\cprime{$'$} \def\cprime{$'$} \def\cprime{$'$} \def\cprime{$'$}
  \def\ocirc#1{\ifmmode\setbox0=\hbox{$#1$}\dimen0=\ht0 \advance\dimen0
  by1pt\rlap{\hbox to\wd0{\hss\raise\dimen0
  \hbox{\hskip.2em$\scriptscriptstyle\circ$}\hss}}#1\else {\accent"17 #1}\fi}
  \def\ocirc#1{\ifmmode\setbox0=\hbox{$#1$}\dimen0=\ht0 \advance\dimen0
  by1pt\rlap{\hbox to\wd0{\hss\raise\dimen0
  \hbox{\hskip.2em$\scriptscriptstyle\circ$}\hss}}#1\else {\accent"17 #1}\fi}

\end{document}